\documentclass[preprint,3p,12pt]{elsarticle}

\usepackage{amssymb, bm}
\usepackage{amsthm, amsmath}

\usepackage{graphicx,color}
\usepackage{caption}
\usepackage{subcaption}
\usepackage{multirow}
\usepackage{mathtools,xparse}

\usepackage{listings}

\makeatletter
\newenvironment{smallermatrix}[1][c]
{\null\,\vcenter\bgroup
  \Let@\restore@math@cr\default@tag
  \baselineskip0pt \lineskip0.4pt \lineskiplimit0pt
  \ialign\bgroup\if#1l\else\hfil\fi$\m@th\scriptstyle##$\if#1r\else\hfil\fi&&\thickspace\hfil
  $\m@th\scriptstyle##$\hfil\crcr
}{%
  \crcr\egroup\egroup\,%
}
\makeatother

\NewDocumentCommand{\ts}{O{c} e{^?_}}{
  \begin{smallermatrix}[#1]
  \mathstrut\IfValueT{#2}{#2} \\
  \mathstrut\IfValueT{#3}{#3} \\
  \mathstrut\IfValueT{#4}{#4}
  \end{smallermatrix}%
}

\newcommand{\Th}{\mathcal{T}_h}

\newcommand{\dTh}{\partial\mathcal{T}_h}

\newcommand{\dV}{\mathrm d \Omega}
\newcommand{\dS}{\mathrm d \Gamma}

\newcommand{\Ohprod}[2]{( #1, #2)_{\Th}}
\newcommand{\pOhprod}[2]{\langle #1, #2\rangle_{\dTh}}

\newcommand{\Kprod}[2]{( #1, #2)_{K}}
\newcommand{\PKprod}[2]{\langle #1, #2\rangle_{\partial K}}

\journal{Arxiv}

\begin{document}

\begin{frontmatter}



\title{Hybridizable discontinuous Galerkin methods for the Monge–Amp\`ere equation}




\author[inst2]{Ngoc Cuong Nguyen}
\author[inst2]{Jaime Peraire}

\affiliation[inst2]{organization={Center for Computational Engineering, Department of Aeronautics and Astronautics, Massachusetts Institute of Technology},
            addressline={77 Massachusetts
Avenue}, 
            city={Cambridge},
            state={MA},
            postcode={02139}, 
            country={USA}}
            

\begin{abstract}
We introduce two hybridizable discontinuous Galerkin (HDG) methods for numerically solving the Monge–Amp\`ere equation.  The first HDG method is devised to solve the nonlinear elliptic Monge–Amp\`ere equation by using  Newton's method. The second HDG method is devised to solve a sequence of the Poisson equation until convergence to a fixed-point solution of the Monge–Amp\`ere equation is reached. Numerical examples are presented to demonstrate the convergence and accuracy of the HDG methods. Furthermore, the HDG methods are applied to  $r$-adaptive mesh generation by redistributing a given scalar density function via the optimal transport theory. This $r$-adaptivity methodology leads to the Monge–Amp\`ere equation with a nonlinear Neumann boundary condition arising from the optimal transport of the density function to conform the resulting high-order mesh to the boundary. Hence, we extend the HDG methods to treat the nonlinear Neumann boundary condition. Numerical experiments are presented to illustrate the generation of $r$-adaptive high-order meshes on planar and curved domains.  
\end{abstract}







\end{frontmatter}


\section{Introduction}
\label{sec:intro}

The Monge-Amp\`ere equation has its root in optimal transport theory and arises from  many areas in science and engineering such as astrophysics, differential geometry, geostrophic fluid dynamics, image processing, mesh generation, optimal transportation, statistical inference, and stochastic control; see \cite{Feng2013} and references therein. The equation belongs to a class of fully nonlinear second-order elliptic partial differential equations (PDEs). Due to the wide range of above-mentioned applications, the Monge-Amp\`ere equation has attracted significant attention from mathematicians and scientists \cite{Brenier1991,Benamou2000,Benamou2010,Benamou2014,Caffarelli1990,Dean2006,Frisch2002,Froese2011,Oliker2003,Prins2014,Trudinger2008}. Significant progress has been made in the development of numerical methods for solving the Monge-Amp\`ere equation \cite{Awanou2015,Benamou2010,Benamou2014,Bohmer2008,Brenner2011,Brenner2012,Caboussat2013,Dean2006,Feng2014,Feng2009,Feng2013,Feng2014a,Feng2017,Feng2018,Froese2011,Froese2012,Liu2019,Lakkis2013}. A number of finite element methods have been proposed for the Monge-Amp\`ere equation. In \cite{Dean2006}, Dean and Glowinski presented an augmented Lagrange multiplier method and a least squares method for the
Monge-Amp\`ere equation by treating the nonlinear equations as a constraint and using finite elements. B\"ohmer \cite{Bohmer2008} introduced a projection method using $C^1$ finite element functions for a class of fully nonlinear second order elliptic PDEs and analyzed the method using consistency and stability arguments. In \cite{Brenner2011}, Brenner et al. proposed $C^0$ finite element methods and discontinuous Galerkin (DG) methods for the Monge–Amp\`ere equation, which were extended to the three dimensional Monge-Amp\`ere equation  \cite{Brenner2012}. In \cite{Benamou2010}, Benamou et al. proposed a fixed-point method that only requires the solution of a sequence of Poisson equations for the two-dimensional Monge–Amp\`ere equation. In \cite{Lakkis2013}, a linearization-then-descretization approach consists in applying  Newton’s method to a nonlinear elliptic PDE to produce a sequence of linear nonvariational elliptic PDEs that can be dealt with using nonvariational finite element method. In \cite{Feng2014a}, Feng and Lewis developed mixed interior penalty DG methods for fully nonlinear second order elliptic and parabolic equations. Recently, a finite element/operator-splitting method is introduced for the Monge–Amp\`ere equation \cite{Glowinski2019,Liu2019} by using an equivalent divergence formulation of the Monge–Amp\`ere equation through the cofactor matrix of the Hessian of the solution.

In this paper, two hybridizable discontinuous Galerkin (HDG) methods are considered for the numerical solution of the Monge–Amp\`ere equation in two dimensions. The first HDG method is devised to solve the Monge–Amp\`ere equation by using  Newton's method. The second HDG method is devised to solve a sequence of the Poisson equation until convergence to a fixed-point solution of the Monge–Amp\`ere equation is reached. Numerical examples are presented to compare the convergence and accuracy of the HDG methods. It is found out that while the two methods yield similar orders of accuracy for the numerical solution, the Newton-HDG method requires considerably less number of linear solves than the fixed-point HDG method. However, an advantage of the fixed-point HDG method is that the discretization of the Poisson equation results in linear systems that can be solved more efficiently than those resulting from the Newton-HDG method. 

The HDG methods have some unique features which distinguish themselve from other finite element methods for the Monge–Amp\`ere equation. First, the global degrees of freedom are those of the approximate trace of the scalar variable defined on element faces. This translates to computational efficiency for the solution of nonlinear and linear systems arising from the HDG discretization of the Monge–Amp\`ere equation. Second, the approximate gradient and Hessian converge with the same order as the approximate scalar variable. For most other finite element methods, the convergence rates of the approximate gradient and Hessian are lower than that of the approximate scalar variable.


There has been considerable interest in $r$-adaptive mesh generation by the optimal transport theory via solving the Monge–Amp\`ere equation \cite{Delzanno2008,Budd2009a,Budd2015,Browne2014,Chacon2011,Weller2016,McRae2018,Sulman2011,Sulman2021}. In $r$-adaptivity, mesh points are neither created nor destroyed, data structures do not need to be modified in-place, and complicated load-balancing is not necessary \cite{Aparicio-Estrems2023}. Furthermore, the $r$-adaptivity approach based on the Monge–Amp\`ere equation has the ability to avoid mesh entanglement and sharp changes in mesh resolution \cite{Budd2009a,Budd2015}. This $r$-adaptivity approach has its root from the optimal transport theory  since it seeks to minimize a deformation functional subject to equidistributing a given scalar monitor function which controls the local density of mesh points. If the functional is defined as the $L_2$ norm of the grid deformation, the optimal transport theory results in a mesh mapping as the gradient of the solution of the Monge–Amp\`ere equation \cite{Delzanno2008,Chacon2011}. To determine the mesh mapping we need to impose a boundary condition that mesh mapping conforms to the boundary of the physical domain. This boundary condition can be characterized as a Hamilton–Jacobi equation on the boundary \cite{Benamou2014} or a nonlinear Neumann boundary condition \cite{Browne2014}. We extend the HDG methods to the Monge–Amp\`ere equation with nonlinear boundary conditions. The methods are used to generate $r$-adaptive high-order meshes based on an equidistribution of a density function via the optimal transport theory. 


The paper is organized as follows. In Section 2, we introduce HDG methods for the Monge-Amp\`ere equation and present numerical results to demonstrate their performance. In Section 3, we extend these HDG methods to $r$-adaptivity based on the optimal transport theory and present numerical experiments to illustrate the generation of $r$-adaptive high-order meshes. Finally, in Section 4, we make a number of concluding remarks on the results as well as future work.

\section{The hybridizable discontinuous Galerkin methods}

HDG methods were first introduced in \cite{CockburnGopalakrishnanLazarov09HDG}  for
elliptic problems and subsequently extended to a wide variety of PDEs: linear convection-diffusion problems~\cite{CockburnDongGuzmanRestelliSacco,Nguyen2009d}, nonlinear convection-diffusion problems~\cite{CockburnDGRS09CDR,Nguyen2009c,Ueckermann2010}, Stokes problems~\cite{CockburnGopalakrishnan09HDGStokes,CockburnGopalakrishnanNguyenPeraireSayasHDGStokes11,Cockburn2010e,Nguyen2010}, incompressible flows equations~\cite{Ahnert2014a,Nguyen2011h,NguyenPeraireCockburn10HDG,NguyenPeraireCockburn10AIAAINS,Rhebergen2012a,Rhebergen2018,Ueckermann2016},  compressible flows \cite{Ciuca2020,Fernandez2017a,Franciolini2020,Moro2011a,Nguyen2012,Nguyen2015c,Vila-Perez2021,Vila-Perez2022}, Maxwell's equations \cite{Nguyen2011j,SANCHEZ2022114969,Li2014}, 
linear elasticity \cite{SoonCockburnStolarski09LE,Cockburn2013,fu2015analysis,qiu2018hdg,Sanchez2021}, and nonlinear elasticity \cite{Cockburn2019,Fernandez2018a,kabaria2015hybridizable,Terrana2019a}. To the best of our knowledge, however, HDG methods have not been considered for solving the Monge–Amp\`ere equation prior to this work. In this section, we describe HDG methods for numerically solving the Monge–Amp\`ere equation. The proposed HDG methods simultaneously compute approximations to the scalar variable, the gradient variables, and Hessian variables. The HDG methods are computationally efficient owing to a hybridization procedure that eliminates the degrees of freedom of those approximate variables to obtain global systems for the degrees of freedom of an approximate trace defined on the element faces.


\subsection{Governing equations and approximation spaces}

We consider the  Monge-Amp\`ere equation with a Dirichlet boundary condition
\begin{equation}
\begin{array}{rcll}
\det (D^2 u) & = & f, \quad & \mbox{in }  \Omega, \\
u & = & g, \quad & \mbox{on } \partial \Omega .
\end{array}
\label{model_problem}
\end{equation}
Here $\Omega \in \mathbb{R}^2$ is the physical domain with Lipschitz boundary $\partial \Omega$, the source term $f$ is {a} strictly positive function, and $g$ is a smooth function. Note that $D^2 u$ is the Hessian matrix of the exact solution $u$, and $\det (D^2 u) = u_{xx} u_{yy} - u_{xy}^2$ is its determinant. The solution $u$ must be convex in order for the equation to be elliptic. Without the convexity constraint, the equation does not have a unique solution. In two dimensions \cite{Benamou2010}, the Monge-Amp\`ere equation (\ref{model_problem}) can be rewritten as
\begin{equation}
\Delta u   =  \sqrt{u_{xx}^2 + u_{yy}^2 + 2 u_{xy}^2 + 2 f} \qquad  \mbox{in }  \Omega .
\label{model_problem2}
\end{equation}
We introduce  $\bm{q} = \nabla u$ and $\bm{H} = \nabla \bm q$, and rewrite the above equation as a first-order system of equations
\begin{equation}
\begin{array}{rcll}
\bm{H} - \nabla \bm q & = & 0, \quad & \mbox{in } \Omega, \\
\bm{q} - \nabla u & = & 0, \quad & \mbox{in } \Omega, \\
s(\bm H, f) - \nabla \cdot \bm q & = & 0 , \quad & \mbox{in }  \Omega, \\
u - g & = & 0, \quad & \mbox{on } \partial \Omega ,
\end{array}
\label{model_problem3}
\end{equation}
where $s(\bm H, f) = \sqrt{H_{11}^2 + H_{22}^2 + H_{12}^2 + H_{21}^2 + 2 f}$. While we consider Dirichlet boundary condition in this section, we will treat Neumann boundary condition in the next section.

We denote by $\mathcal{T}_h$ a collection of disjoint regular elements $K$ that partition $\Omega$ and set $\partial\mathcal{T}_h := \{\partial K:\;K\in\mathcal{T}_h\}$.
For an element $K$ of the collection $\mathcal{T}_h$, $e = \partial K \cap \partial \Omega$ is the boundary face if the {$(d-1)$-Lebesgue} measure of $e$ is nonzero. For two elements $K^+$ and $K^-$ of the collection $\mathcal{T}_h$, $e = \partial K^+ \cap \partial K^-$ is the interior face between $K^+$ and $K^-$ if the {$(d-1)$-Lebesgue} measure of $e$ is nonzero. Let $\mathcal{E}_h^i$ and $\mathcal{E}_h^{\partial}$ denote the set of interior and boundary faces, respectively. We denote by $\mathcal{E}_h$ the union of $\mathcal{E}_h^{i}$ and $\mathcal{E}_h^{\partial}$. Let $\bm{n}^+$ and $\bm{n}^-$ be the outward unit normal of $\partial K^+$ and $\partial K^-$, respectively. Let $\mathcal{P}^p(D)$ denote the set of polynomials of degree at most $p$ on a domain $D$ and let $L^2(D)$ be the space of square integrable functions on $D$. We introduce discontinuous finite element spaces
\begin{equation*}
\begin{split}
\bm{W}_{h}^p &  = \{\bm{v} \in \bm{L}^2(\Omega) \ : \ \bm{v}|_K \in (\mathcal{P}^p(K))^{2 \times 2} \quad \forall K \in \mathcal{T}_h \}, \\
\bm{V}_{h}^p &  = \{\bm{v} \in \bm{L}^2(\Omega) \ : \ \bm{v}|_K \in (\mathcal{P}^p(K))^2 \quad \forall K \in \mathcal{T}_h \}, \\
U_{h}^p  &  = \{w \in L^2(\Omega) \ : \ w|_K \in \mathcal{P}^p(K) \quad \forall K \in \mathcal{T}_h \}, \\
M_{h}^p  & = \{\mu \in L^2(\mathcal{E}_h) \ : \ \mu|_e \in \mathcal{P}^p(e), \quad \ \forall e \in \mathcal{E}_h \} .
\end{split}
\end{equation*}
Note that $M_h^p$ consists of functions which are continuous inside the faces (or edges)
$e \in \mathcal{E}_h$ and discontinuous at their borders. 


For functions $u$ and $v$ in $L^2(D)$, we denote $(u,v)_D = \int_{D} u v$ if $D$ is a domain in $\mathbb{R}^d$ and $\left\langle u,v\right\rangle_D = \int_{D} u v$ if $D$ is a domain in $\mathbb{R}^{d-1}$. The inner produces associated with the above approximation spaces are defined as follows
\begin{subequations}\label{eq:inner_products}
\begin{alignat}{5}
\Kprod{{u}}{{v}} &:= \int_{K}{u} {v}\,\dV , &\quad \Ohprod{{u}}{{v}} &:= \sum_{K\in \Th}\Kprod{{u}}{{v}},
   &\quad \forall {u}, {v}\in L^2(\Omega), \nonumber \\
  \Kprod{\bm{u}}{\bm{v}} &:= \int_{K}\bm{u}\cdot\bm{v}\,\dV , &\quad \Ohprod{\bm{u}}{\bm{v}} &:= \sum_{K\in \Th}\Kprod{\bm{u}}{\bm{v}},
   &\quad \forall \bm{u},\bm{v}\in L^2(\Omega)^d, \nonumber \\
  \Kprod{\bm{G}}{\bm{H}} &:= \int_{K} \bm{G}:\bm{H} \,\dV , &\quad \Ohprod{\bm{G}}{\bm{H}} &:= \sum_{K\in \Th} \Kprod{\bm{G}}{\bm{H}},
  &\quad \forall \bm{G},\bm{H}\in L^2(\Omega)^{d\times d},  \nonumber  \\
\PKprod{{\mu}}{{\eta}} &:= \int_{\partial K} {\mu} {\eta} \,\dS , &\quad \pOhprod{{\mu}}{{\eta}} &:= \sum_{K\in \Th} \PKprod{{\mu}}{{\eta}}, &\quad  \forall {\mu},{\eta}\in L^2(\dTh).  \nonumber 
\end{alignat}
\end{subequations}
We are ready to describe the HDG methods for the Monge-Amp\`ere problem~\eqref{model_problem3}. 

\subsection{The Newton-HDG formulation}

To numerically solve \eqref{model_problem3} with the HDG method, we find $(\bm{H}_h, \bm{q}_h,u_h,\widehat{u}_h) \in \bm{W}_{h}^p \times \bm{V}_{h}^p  \times U_{h}^p \times M_h^p$ such that
\begin{equation}
\label{eq8}
\begin{array}{rcl}
 \left(\bm{H}_h, \bm{G}\right)_{\mathcal{T}_h} +  \left(\bm q_h,  \nabla \cdot \bm{G}\right)_{\mathcal{T}_h} -  \left\langle \widehat{\bm q}_h, \bm{G} \cdot \bm{n} \right\rangle_{\partial \mathcal{T}_h} & = & 0,  \\
 \left(\bm{q}_h, \bm{v}\right)_{\mathcal{T}_h} +  \left(u_h,  \nabla \cdot \bm{v}\right)_{\mathcal{T}_h} -  \left\langle \widehat{u}_h, \bm{v} \cdot \bm{n} \right\rangle_{\partial \mathcal{T}_h} & = & 0,  \\
\left(\bm{q}_h, \nabla w\right)_{\mathcal{T}_h}  -  \left\langle \widehat{\bm{q}}_h \cdot \bm{n}, w \right\rangle_{\partial \mathcal{T}_h} + (s(\bm H_h, f),w)_{\mathcal{T}_h} & = & 0, \\
\left\langle \widehat{\bm{q}}_h  \cdot \bm{n} , \mu \right\rangle_{\partial \mathcal{T}_h \backslash \partial \Omega} + \left\langle \widehat{u}_h  - g, \mu \right\rangle_{\partial \Omega} & =  & 0,
\end{array}
\end{equation}
for all $(\bm G, \bm{v}, w, \mu) \in  \bm{W}_h^p \times 
 \bm{V}_h^p \times U_h^p \times M_h^p$, where 
\begin{equation}
\widehat{\bm{q}}_h = {\bm{q}_h} -
\tau (u_h - \widehat{u}_h) \bm{n}, \quad \mbox{on } \mathcal{E}_h.
\label{fluxdef}
\end{equation}
Here $\tau$ is the stabilization parameter which is set to 1 for all the numerical examples presented in this paper. 


Newton's method is used to solve the nonlinear system~\eqref{eq8}. The procedure evaluates successive approximations $(\bm{H}_h^l, \bm{q}_h^l,u_h^l,\widehat{u}_h^l)$ starting from an initial guess $(\bm{H}_h^0, \bm{q}_h^0,u_h^0,\widehat{u}_h^0)$. For each Newton step $l$, the system of equations \eqref{eq8} is linearized with respect to the Newton increments $\left(\delta \bm{H}_h^l, \delta  \bm{q}_h^l, \delta u_h^l, \delta 
 \widehat{u}_h^l \right) \in \bm{W}_{h}^p \times \bm{V}_{h}^p  \times U_{h}^p \times M_h^p$ that satisfy 
\begin{subequations}\label{eq:HDG_increments}
\begin{alignat}{1}
\left(\delta \bm{H}_h^l, \bm{G}\right)_{\mathcal{T}_h} +  \left(\delta \bm q_h^l,  \nabla \cdot \bm{G}\right)_{\mathcal{T}_h} -  \left\langle \delta \widehat{\bm q}_h^l, \bm{G} \cdot \bm{n} \right\rangle_{\partial \mathcal{T}_h}  & =  r_1(\bm G) , \label{subeq:increm1} \\
  \left(\delta \bm{q}_h^l, \bm{v}\right)_{\mathcal{T}_h} +  \left(\delta u_h^l,  \nabla \cdot \bm{v}\right)_{\mathcal{T}_h} -  \left\langle \delta \widehat{u}_h^l, \bm{v} \cdot \bm{n} \right\rangle_{\partial \mathcal{T}_h}  & = r_2(\bm{v}) , \label{subeq:increm2} \\
  (\partial s_{\bm H}(\bm H_h^l, f) \delta \bm H_h^l ,w)_{\mathcal{T}_h} + \left(\delta \bm{q}_h^l, \nabla w\right)_{\mathcal{T}_h}  -  \left\langle \delta \widehat{\bm{q}}_h^l \cdot \bm{n}, w \right\rangle_{\partial \mathcal{T}_h}  & = r_3(w) , \label{subeq:increm3} \\
  \left\langle \delta \widehat{\bm{q}}_h^l  \cdot \bm{n} , \mu \right\rangle_{\partial \mathcal{T}_h \backslash \partial \Omega}  + \left\langle \delta \widehat{u}_h^l , \mu \right\rangle_{\partial \Omega}  & = r_4({\mu}), \label{subeq:increm4}
\end{alignat}
\end{subequations}
for all $(\bm G, \bm{v}, w, \mu) \in  \bm{W}_h^p \times 
 \bm{V}_h^p \times U_h^p \times M_h^p$, the right-hand side residuals are given by
\begin{subequations}\label{eq:HDG_residuals}
\begin{alignat}{1}
  r_1(\bm G)     & = - \left(\bm{H}_h^l, \bm{G}\right)_{\mathcal{T}_h} -  \left(\bm q^l_h,  \nabla \cdot \bm{G}\right)_{\mathcal{T}_h} +  \left\langle \widehat{\bm q}_h^l, \bm{G} \cdot \bm{n} \right\rangle_{\partial \mathcal{T}_h} \\
  r_2(\bm{v})  & =- \left(\bm{q}_h^l, \bm{v}\right)_{\mathcal{T}_h} -  \left(u_h^l,  \nabla \cdot \bm{v}\right)_{\mathcal{T}_h} +  \left\langle \widehat{u}_h^l, \bm{v} \cdot \bm{n} \right\rangle_{\partial \mathcal{T}_h} \label{subeq:resid2} \\
  r_3(w)& = - \left(\bm{q}_h^l, \nabla w\right)_{\mathcal{T}_h}  +  \left\langle \widehat{\bm{q}}_h^l \cdot \bm{n}, w \right\rangle_{\partial \mathcal{T}_h} - (s(\bm H_h^l, f),w)_{\mathcal{T}_h}   \label{subeq:resid3} \\
r_4(\mu) & = -\left\langle \widehat{\bm{q}}_h^l  \cdot \bm{n} , \mu \right\rangle_{\partial \mathcal{T}_h \backslash \partial \Omega} - \left\langle \widehat{u}_h^l  - g, \mu \right\rangle_{\partial \Omega} .  \label{subeq:resid4}
\end{alignat}
\end{subequations}
Here $\partial s_{\bm H}(\bm H, f)$ denotes the partial derivatives of $s(\bm H, f)$ with respect to $\bm H$. After solving \eqref{eq:HDG_increments}, the numerical approximations are then updated
\begin{equation}\label{eq:update_increm}
  \left( \bm{H}_h^{l+1}, \bm{q}_h^{l+1},u_h^{l+1},\widehat{u}_h^{l+1} \right)  := \left( \bm{H}_h^l, \bm{q}_h^l,u_h^l,\widehat{u}_h^l \right) + \alpha \left( \delta \bm{H}_h^l, \delta  \bm{q}_h^l, \delta u_h^l, \delta 
 \widehat{u}_h^l  \right),
\end{equation}
where the coefficient $\alpha$ is determined by a line-search algorithm in order to optimally decrease the residual. This process is repeated until the residual norm is smaller than a given tolerance, typically $10^{-8}$. 

At each step of the Newton method, the linearization \eqref{eq:HDG_increments} gives the following matrix system to be solved
\begin{equation}\label{eq:matrix_system}
  \begin{pmatrix} \mathbb{A}^l & \mathbb{B}^l \\ \mathbb{C}^l & \mathbb{D}^l \end{pmatrix} \begin{pmatrix} \delta U^l  \\ \delta \widehat{U}^l \end{pmatrix}
    = \begin{pmatrix} R_{123}^l  \\ R_4^l \end{pmatrix} ,
\end{equation}
where $\delta U^l$ and $\delta \widehat{U}^l$ are the vectors of degrees of freedom of $(\delta \bm{H}_h^l, \delta  \bm{q}_h^l, \delta u_h^l)$ and $\delta 
 \widehat{u}_h^l$, respectively. The system \eqref{eq:matrix_system} is first solved for the traces only $\delta \widehat{U}^l$ as
\begin{equation}\label{eq:reduc_matrix_system}
  \mathbb{K}^l \delta \widehat{U}^l = R^l,
\end{equation}
where $\mathbb{K}^l$ is the Schur complement of the block $\mathbb{A}^l$ and $R^l$ is the reduced residual
\begin{equation}\label{eq:SchurCompl}
  \mathbb{K}^l = \mathbb{D}^l - \mathbb{C}^l \left( \mathbb{A}^l \right)^{-1} \mathbb{B}^l , \ \ \ \ \ \
  R^l = R_4^l - \mathbb{C}^l \left( \mathbb{A}^l \right)^{-1} R_{123}^l .
\end{equation}
The reduced system \eqref{eq:reduc_matrix_system} involves fewer degrees of freedom than the full system \eqref{eq:matrix_system}. Moreover due to the discontinuous nature of the approximate solution, the matrix $\mathbb{A}^l$ and its inverse are block diagonal, and can be computed elementwise. Once $\delta \widehat{U}^l$ is known, the other unknowns $\delta U^l$ are then retrieved element-wise. Therefore, the full system \eqref{eq:matrix_system} is never explicitly built, and the reduced matrix $\mathbb{K}^l$ is built directly in an elementwise fashion, thus reducing the memory storage.

\subsection{The fixed point-HDG method}

In addition to using Newton's method, we employ the fixed-point method to solve  the nonlinear system~\eqref{eq8}. We note that the nonlinear source term $s$ in \eqref{eq8} depends on $\bm H_h$. Hence, starting from an initial guess $\bm H^0_h$  we repeatedly find $(\bm{q}_h^l,u_h^l,\widehat{u}_h^l) \in \bm{V}_{h}^p  \times U_{h}^p \times M_h^p$ such that
\begin{equation}
\label{fpHDG}
\begin{array}{rcl}
  \left(\bm{q}_h^l, \bm{v}\right)_{\mathcal{T}_h} +  \left(u_h^l,  \nabla \cdot \bm{v}\right)_{\mathcal{T}_h} -  \left\langle \widehat{u}_h^l, \bm{v} \cdot \bm{n} \right\rangle_{\partial \mathcal{T}_h} & = & 0,  \\
\left(\bm{q}_h^l, \nabla w\right)_{\mathcal{T}_h}  -  \left\langle \widehat{\bm{q}}_h^l \cdot \bm{n}, w \right\rangle_{\partial \mathcal{T}_h}  & = & - (s(\bm H_h^{l-1}, f),w)_{\mathcal{T}_h}, \\
\left\langle \widehat{\bm{q}}_h^l  \cdot \bm{n} , \mu \right\rangle_{\partial \mathcal{T}_h \backslash \partial \Omega} + \left\langle \widehat{u}_h^l , \mu \right\rangle_{\partial \Omega} & =  &  \left\langle  g, \mu \right\rangle_{\partial \Omega},
\end{array}
\end{equation}
for all $(\bm{v}, w, \mu) \in  \bm{V}_h^p \times U_h^p \times M_h^p$, and then compute $\bm{H}_h^l\in \bm{W}_{h}^p$ such that
\begin{equation}
\label{eqH}
 \left(\bm{H}_h^l, \bm{G}\right)_{\mathcal{T}_h} =  - \left(\bm q_h^l,  \nabla \cdot \bm{G}\right)_{\mathcal{T}_h} +  \left\langle \widehat{\bm q}_h^l, \bm{G} \cdot \bm{n} \right\rangle_{\partial \mathcal{T}_h}, \quad \forall \bm G \in \bm{W}_h^p 
\end{equation}
until $\|\bm H^l_h - \bm H^{l-1}_h\|$ is less than a specified tolerance, typically $10^{-6}$. 

We note that the weak formulation (\ref{fpHDG}) is nothing but the HDG method for the Poisson equation $-\Delta u = -s$ in $\Omega$ with $u = g$ on $\partial \Omega$. Applying the solution strategy described earlier to the system (\ref{fpHDG}), we obtain the following global system 
\begin{equation}
\label{eqfp:reduc_matrix_system}
  \mathbb{G} \  \widehat{U}^l = F^l,
\end{equation}
for the degrees of freedom of $\widehat{u}_h^l$. As the global matrix $\mathbb{G}$ is unchanged, it can be computed once prior to carrying out the fixed-point algorithm. However, the right hand side vector $F^l$ has to be computed at every fixed-point iteration because the source term depends on the previous iterate. Once $ \widehat{U}^l$ is obtained by solving the linear system (\ref{eqfp:reduc_matrix_system}), the degrees of freedom of both $\bm q_h^l$ and $u_h^l$ are then retrieved element-wise. Finally, the degrees of freedom of $\bm H^l_h$ are also computed by solving (\ref{eqH}) element-wise. 

We see that the fixed-point HDG method for the Monge-Amp\`ere equation means to solving the Poisson problem with a sequence of right-hand side vectors. Hence, the fixed-point HDG method is much easier to implement the Newton-HDG method. Furthermore, its global matrix is computed only once, whereas that of the Newton-HDG method has to be computed for each Newton step. However, the Newton-HDG method can converge considerably faster than the fixed-point HDG method.

\subsection{Numerical examples}

In this section we provide examples to demonstrate the HDG methods described above. In the first example we compare the convergence and accuracy of the methods for several polynomial degrees. In the second example we focus on a problem with a nearly singular solution. These examples demonstrate the relative advantages and disadvantages of the HDG methods. For the first problem, the Newton-HDG method has robust convergence in the nonlinear iteration, while the fixed-point HDG method requires more iterations. The second example illustrates the influence of regularity on the performance of the methods. For these examples, we choose $\mathcal{T}_h$ as a uniform triangulation of $\Omega$ with mesh size  $h = 1/n$ and use $u^0 = (x^2+y^2)/2$ as the initial guess.

\

\noindent
\textbf{Example 1.} We consider the Monge-Amp\`ere equation in $\Omega = (0,1)^2$ with $f = (1 + x^2 +y^2) e^{x^2 + y^2}$ and $g = e^{0.5 (x^2 + y^2)}$, which yields the exact solution $u = e^{0.5 (x^2 + y^2)}$. We use both triangular and quadrilateral meshes with resolutions of $n = 4, 8, 16, 32$, $64$, and polynomial degrees of $p = 1, 2$, $3$.

We report the $L^2$ errors  and orders of convergence of the computed solutions on the triangular meshes in Table \ref{ex1tab1} and the quadrilateral meshes in Table \ref{ex1tab2}. We see that both the Newton-HDG method and the fixed-point HDG method have the same errors and convergence rates. These results are expected because the two methods should converge to the same numerical solution. We observe interestingly that the convergence rates of $u_h$, $\bm q_h$, and $\bm H_h$ are $O(h^p)$ on  triangular and quadrilateral meshes.  The convergence rate of $O(h^p)$ for the computed Hessian $\bm H_h$ is a good news for the HDG methods since the Hessian $H$ are the second-order partial derivatives of the scalar variable $u$. However, the convergence rate of $O(h^p)$ for both $u_h$ and $\bm q_h$ are suboptimal. The suboptimal convergence of $u_h$ and $\bm q_h$ for the Monge-Amp\`ere equation can be attributed to the fact that the source term $s$ depends on the Hessian $\bm H$. As the computed Hessian $\bm H_h$ converges with order $O(h^p)$, the $L^2$ projection of the source term $s(\bm H_h, f)$ also converges with order $O(h^p)$. In contrast, for linear second-order elliptic problems, the $L^2$ projection of the source term converges with order $O(h^{p+1})$, the HDG method can yield optimal convergence rate of $O(h^{p+1})$ for both $u_h$ and $\bm q_h$. Lastly, we show the number of iterations required to reach convergence in Table \ref{ex1tab3}. As expected, the Newton-HDG method requires considerably less iterations (about six times less) than the fixed-point HDG method. Furthermore, the number of iterations appears consistent for all values of $n$ and $p$ on both triangular and quadrilateral meshes.

\begin{table}[ht]
  \begin{center}
\scalebox{0.8}{%
    $\begin{array}{|c|c||c c  | c c|  c c ||  c c| c c| c c|}
    \hline    
    & & \multicolumn{6}{|c||}{\mbox{Newton-HDG method}} & \multicolumn{6}{|c|}{\mbox{Fixed point-HDG method}} \\ 
   \hline 
  \mbox{degree} & \mbox{mesh} & \multicolumn{2}{|c}{\|u-u_h\|} &  \multicolumn{2}{|c}{\|\bm{q}-\bm{q}_h\|} & \multicolumn{2}{|c||}{\|\bm{H}-\bm{H}_h\|} & \multicolumn{2}{|c|}{\|u-u_h\|} &  \multicolumn{2}{|c}{\|\bm{q}-\bm{q}_h\|} & \multicolumn{2}{|c|}{\|\bm{H}-\bm{H}_h\|}  \\
    p & n & \mbox{error}  & \mbox{order} & \mbox{error}  & \mbox{order} & \mbox{error}  & \mbox{order} & \mbox{error}  & \mbox{order} & \mbox{error}  & \mbox{order} & \mbox{error}  & \mbox{order}  \\
  \hline    
  & 4  &  1.37\mbox{e-}3  &  --  &  1.12\mbox{e-}2  &  --  &  5.75\mbox{e-}1  &  --  &  1.37\mbox{e-}3  &  --  &  1.12\mbox{e-}2  &  --  &  5.75\mbox{e-}1  &  --  \\  
 & 8  &  1.13\mbox{e-}4  &  1.96  &  2.73\mbox{e-}3  &  2.04  &  2.95\mbox{e-}1  &   0.97  &  1.13\mbox{e-}4  &  1.96  &  2.73\mbox{e-}3  &  2.04  &  2.95\mbox{e-}1  &   0.97  \\  
1 & 16  &  1.75\mbox{e-}4  &  1.80  &  1.16\mbox{e-}3  &  1.23  &  1.49\mbox{e-}1  &   0.98  &  1.75\mbox{e-}4  &  1.80  &  1.16\mbox{e-}3  &  1.23  &  1.49\mbox{e-}1  &   0.98  \\  
 & 32  &  1.37\mbox{e-}4  &  1.25  &  7.30\mbox{e-}4  &   0.67  &  7.49\mbox{e-}2  &   0.99  &  1.37\mbox{e-}4  &  1.25  &  7.30\mbox{e-}4  &   0.67  &  7.49\mbox{e-}2  &   0.99  \\  
 & 64  &  8.13\mbox{e-}5  &   0.94  &  4.21\mbox{e-}4  &   0.79  &  3.76\mbox{e-}2  &  1.00  &  8.13\mbox{e-}5  &   0.94  &  4.21\mbox{e-}4  &   0.79  &  3.76\mbox{e-}2  &  1.00  \\  
\hline
 & 4  &  1.29\mbox{e-}4  &  --  &  1.01\mbox{e-}3  &  --  &  5.82\mbox{e-}2  &  --  &  1.29\mbox{e-}4  &  --  &  1.01\mbox{e-}3  &  --  &  5.82\mbox{e-}2  &  --  \\  
 & 8  &  4.28\mbox{e-}5  &  2.51  &  2.44\mbox{e-}4  &  2.05  &  1.49\mbox{e-}2  &  1.96  &  4.28\mbox{e-}5  &  2.51  &  2.44\mbox{e-}4  &  2.05  &  1.49\mbox{e-}2  &  1.96  \\  
2 & 16  &  1.17\mbox{e-}5  &  2.21  &  6.26\mbox{e-}5  &  1.96  &  3.78\mbox{e-}3  &  1.98  &  1.17\mbox{e-}5  &  2.21  &  6.26\mbox{e-}5  &  1.96  &  3.78\mbox{e-}3  &  1.98  \\  
 & 32  &  3.03\mbox{e-}6  &  2.05  &  1.60\mbox{e-}5  &  1.97  &  9.49\mbox{e-}4  &  1.99  &  3.03\mbox{e-}6  &  2.05  &  1.60\mbox{e-}5  &  1.97  &  9.49\mbox{e-}4  &  1.99  \\  
 & 64  &  7.67\mbox{e-}7  &  2.01  &  4.03\mbox{e-}6  &  1.98  &  2.38\mbox{e-}4  &  2.00  &  7.67\mbox{e-}7  &  2.01  &  4.03\mbox{e-}6  &  1.98  &  2.38\mbox{e-}4  &  2.00  \\  
\hline
 & 4  &  1.62\mbox{e-}6  &  --  &  4.22\mbox{e-}5  &  --  &  4.37\mbox{e-}3  &  --  &  1.62\mbox{e-}6  &  --  &  4.22\mbox{e-}5  &  --  &  4.37\mbox{e-}3  &  --  \\  
 & 8  &  3.13\mbox{e-}7  &  3.82  &  3.23\mbox{e-}6  &  3.71  &  5.60\mbox{e-}4  &  2.96  &  3.13\mbox{e-}7  &  3.82  &  3.23\mbox{e-}6  &  3.71  &  5.60\mbox{e-}4  &  2.96  \\  
3 & 16  &  5.35\mbox{e-}8  &  3.60  &  3.34\mbox{e-}7  &  3.27  &  7.04\mbox{e-}5  &  2.99  &  5.35\mbox{e-}8  &  3.60  &  3.34\mbox{e-}7  &  3.27  &  7.04\mbox{e-}5  &  2.99  \\  
 & 32  &  7.59\mbox{e-}9  &  3.23  &  4.17\mbox{e-}8  &  3.00  &  8.80\mbox{e-}6  &  3.00  &  7.59\mbox{e-}9  &  3.23  &  4.17\mbox{e-}8  &  3.00  &  8.80\mbox{e-}6  &  3.00  \\  
 & 64  &  9.9\mbox{e-}10  &  3.05  &  5.31\mbox{e-}9  &  2.98  &  1.10\mbox{e-}6  &  3.00  &  9.9\mbox{e-}10  &  3.05  &  5.31\mbox{e-}9  &  2.98  &  1.10\mbox{e-}6  &  3.00  \\  
\hline
     \end{array} $
}
\caption{Error and convergence rate for the HDG methods on the triangular meshes for the Example 1.}
\label{ex1tab1}
\end{center}{$\phantom{|}$}     
\end{table}

\begin{table}[ht]
  \begin{center}
\scalebox{0.8}{%
    $\begin{array}{|c|c||c c  | c c|  c c ||  c c| c c| c c|}
    \hline    
    & & \multicolumn{6}{|c||}{\mbox{Newton-HDG method}} & \multicolumn{6}{|c|}{\mbox{Fixed point-HDG method}} \\ 
   \hline 
  \mbox{degree} & \mbox{mesh} & \multicolumn{2}{|c}{\|u-u_h\|} &  \multicolumn{2}{|c}{\|\bm{q}-\bm{q}_h\|} & \multicolumn{2}{|c||}{\|\bm{H}-\bm{H}_h\|} & \multicolumn{2}{|c|}{\|u-u_h\|} &  \multicolumn{2}{|c}{\|\bm{q}-\bm{q}_h\|} & \multicolumn{2}{|c|}{\|\bm{H}-\bm{H}_h\|}  \\
    p & n & \mbox{error}  & \mbox{order} & \mbox{error}  & \mbox{order} & \mbox{error}  & \mbox{order} & \mbox{error}  & \mbox{order} & \mbox{error}  & \mbox{order} & \mbox{error}  & \mbox{order}  \\
  \hline    
& 4  &  2.20\mbox{e-}3  &  --  &  2.64\mbox{e-}2  &  --  &  5.03\mbox{e-}1  &  --  &  2.20\mbox{e-}3  &  --  &  2.64\mbox{e-}2  &  --  &  5.03\mbox{e-}1  &  --  \\  
 & 8  &  4.68\mbox{e-}4  &  1.73  &  8.85\mbox{e-}3  &  1.58  &  2.92\mbox{e-}1  &   0.78  &  4.68\mbox{e-}4  &  1.73  &  8.85\mbox{e-}3  &  1.58  &  2.92\mbox{e-}1  &   0.78  \\  
1 & 16  &  1.64\mbox{e-}4  &  1.83  &  3.05\mbox{e-}3  &  1.54  &  1.66\mbox{e-}1  &   0.82  &  1.64\mbox{e-}4  &  1.83  &  3.05\mbox{e-}3  &  1.54  &  1.66\mbox{e-}1  &   0.82  \\  
 & 32  &  1.04\mbox{e-}4  &  1.69  &  1.19\mbox{e-}3  &  1.36  &  9.43\mbox{e-}2  &   0.82  &  1.04\mbox{e-}4  &  1.69  &  1.19\mbox{e-}3  &  1.36  &  9.43\mbox{e-}2  &   0.82  \\  
 & 64  &  6.45\mbox{e-}5  &  1.27  &  5.28\mbox{e-}4  &  1.17  &  5.40\mbox{e-}2  &   0.80  &  6.45\mbox{e-}5  &  1.27  &  5.28\mbox{e-}4  &  1.17  &  5.40\mbox{e-}2  &   0.80  \\  
\hline
 & 4  &  1.04\mbox{e-}4  &  --  &  1.62\mbox{e-}3  &  --  &  5.49\mbox{e-}2  &  --  &  1.04\mbox{e-}4  &  --  &  1.62\mbox{e-}3  &  --  &  5.49\mbox{e-}2  &  --  \\  
 & 8  &  3.20\mbox{e-}5  &  2.63  &  2.89\mbox{e-}4  &  2.49  &  1.62\mbox{e-}2  &  1.76  &  3.20\mbox{e-}5  &  2.63  &  2.89\mbox{e-}4  &  2.49  &  1.62\mbox{e-}2  &  1.76  \\  
2 & 16  &  8.56\mbox{e-}6  &  2.48  &  5.70\mbox{e-}5  &  2.34  &  4.8\mbox{e-}3  &  1.76  &  8.56\mbox{e-}6  &  2.48  &  5.70\mbox{e-}5  &  2.34  &  4.8\mbox{e-}3  &  1.76  \\  
 & 32  &  2.20\mbox{e-}6  &  2.24  &  1.26\mbox{e-}5  &  2.18  &  1.42\mbox{e-}3  &  1.75  &  2.20\mbox{e-}6  &  2.24  &  1.26\mbox{e-}5  &  2.18  &  1.42\mbox{e-}3  &  1.75  \\  
 & 64  &  5.55\mbox{e-}7  &  2.08  &  2.97\mbox{e-}6  &  2.08  &  4.24\mbox{e-}4  &  1.75  &  5.55\mbox{e-}7  &  2.08  &  2.97\mbox{e-}6  &  2.08  &  4.24\mbox{e-}4  &  1.75  \\  
\hline
 & 4  &  1.93\mbox{e-}6  &  --  &  5.88\mbox{e-}5  &  --  &  3.24\mbox{e-}3  &  --  &  1.93\mbox{e-}6  &  --  &  5.88\mbox{e-}5  &  --  &  3.24\mbox{e-}3  &  --  \\  
 & 8  &  2.11\mbox{e-}7  &  3.80  &  5.09\mbox{e-}6  &  3.53  &  4.93\mbox{e-}4  &  2.72  &  2.11\mbox{e-}7  &  3.80  &  5.09\mbox{e-}6  &  3.53  &  4.93\mbox{e-}4  &  2.72  \\  
3 & 16  &  2.95\mbox{e-}8  &  3.79  &  4.65\mbox{e-}7  &  3.45  &  7.41\mbox{e-}5  &  2.73  &  2.95\mbox{e-}8  &  3.79  &  4.65\mbox{e-}7  &  3.45  &  7.41\mbox{e-}5  &  2.73  \\  
 & 32  &  4.32\mbox{e-}9  &  3.55  &  4.67\mbox{e-}8  &  3.31  &  1.11\mbox{e-}5  &  2.74  &  4.32\mbox{e-}9  &  3.55  &  4.67\mbox{e-}8  &  3.31  &  1.11\mbox{e-}5  &  2.74  \\  
 & 64  &  6.0\mbox{e-}10  &  3.19  &  5.15\mbox{e-}9  &  3.18  &  1.65\mbox{e-}6  &  2.74  &  6.0\mbox{e-}10  &  3.19  &  5.14\mbox{e-}9  &  3.18  &  1.66\mbox{e-}6  &  2.74  \\  
\hline
     \end{array} $
}
\caption{Error and convergence rate for the HDG methods on the quadrilateral meshes for the Example 1.}
\label{ex1tab2}
\end{center}{$\phantom{|}$}     
\end{table}

\begin{table}[ht]
  \begin{center}
\scalebox{0.8}{%
    $\begin{array}{|c||c c | c c|  c c ||  c c| c c| c c|}
    \hline    
    & \multicolumn{6}{|c||}{\mbox{Newton-HDG method}} & \multicolumn{6}{|c|}{\mbox{Fixed point-HDG method}} \\ 
   \hline 
   \mbox{mesh} & \multicolumn{2}{|c}{p=1} &  \multicolumn{2}{|c}{p=2} & \multicolumn{2}{|c||}{p=3} & \multicolumn{2}{|c|}{p=1} &  \multicolumn{2}{|c}{p=2} & \multicolumn{2}{|c|}{p=3}  \\
     n & \mbox{tri}  & \mbox{quad} & \mbox{tri}  & \mbox{quad} &  \mbox{tri}  &  \mbox{quad} & \mbox{tri}  &  \mbox{quad} & \mbox{tri}  &  \mbox{quad} & \mbox{tri} & \mbox{quad}  \\
  \hline    
4  &  6  &  6  &  6  &  6  &  6  &  6  &  30  &  37  &  38  &  38  &  40  &  42  \\  
 8  &  6  &  6  &  6  &  6  &  6  &  6  &  36  &  40  &  41  &  41  &  41  &  43  \\  
 16  &  7  &  6  &  6  &  6  &  6  &  6  &  40  &  43  &  42  &  43  &  42  &  44  \\  
 32  &  7  &  7  &  6  &  6  &  6  &  6  &  42  &  44  &  43  &  44  &  44  &  45  \\  
 64  &  7  &  7  &  6  &  6  &  6  &  6  &  44  &  44  &  45  &  45  &  45  &  46  \\   
\hline
     \end{array} $
}
\caption{Number of iterations required to reach convergence for the HDG methods for the Example 1. Here ``tri'' columns refer to triangular meshes, while ``quad'' columns refer to quadrilateral meshes. }
\label{ex1tab3}
\end{center}{$\phantom{|}$}     
\end{table}

\

\noindent
\textbf{Example 2.} We consider the Monge-Amp\`ere equation in $\Omega = (0,1)^2$ with $f = R^2/(R^2 - x^2 - y^2)^2$ and $g = - \sqrt{R^2 - x^2 - y^2}$, which yields the exact solution $u = - \sqrt{R^2 - x^2 - y^2}$ for $R > \sqrt{2}$. Note that the regularity of the exact solution decreases as $R$ decreases toward $\sqrt{2}$. Indeed, the gradient and Hessian of the exact solution are singular at the corner $(1,1)$ for $R = \sqrt{2}$.  In this example, we would like to see how the HDG methods perform as the regularity of the solution decreases. Hence, we will report numerical results for $R = 2, R = \sqrt{2} + 0.1,$ and $R = \sqrt{2} + 0.01$. We use triangular meshes with resolutions of $n = 4, 8, 16, 32$, $64, 128, 256$, and polynomial degrees of $p = 1, 2$, $3$.

We report the $L^2$ errors  and orders of convergence of the computed solutions in Table \ref{ex2tab1} for $R=2$, Table \ref{ex2tab2} for $R=\sqrt{2}+0.1$, and Table \ref{ex2tab3} for $R=\sqrt{2}+0.01$. The Newton-HDG method and the fixed-point HDG method have very similar errors and convergence rates. We observe again that $u_h$, $\bm q_h$, and $\bm H_h$ converge with order $O(h^p)$, except for $p=3$ and $R=\sqrt{2}+0.1$ where both $u_h$ and $\bm q_h$ appear to converge like $O(h^{p+1})$. As $R$ decreases, the $L^2$ errors tend to increase because both $\bm q$ and $\bm H$ become less smooth. Hence, we need to increase the mesh resolution in order to observe the expected convergence rates for $R = \sqrt{2}+0.01$. Table \ref{ex2tab4} lists the number of iterations required to reach convergence. The Newton-HDG method requires considerably less iterations than the fixed-point HDG method. The number of iterations for the Newton-HDG method remains consistent, whereas that for the fixed-point HDG method tends to increase as $R$ decreases.

\begin{table}[thbp]
  \begin{center}
\scalebox{0.8}{%
    $\begin{array}{|c|c||c c  | c c|  c c ||  c c| c c| c c|}
    \hline    
    & & \multicolumn{6}{|c||}{\mbox{Newton-HDG method}} & \multicolumn{6}{|c|}{\mbox{Fixed point-HDG method}} \\ 
   \hline 
  \mbox{degree} & \mbox{mesh} & \multicolumn{2}{|c}{\|u-u_h\|} &  \multicolumn{2}{|c}{\|\bm{q}-\bm{q}_h\|} & \multicolumn{2}{|c||}{\|\bm{H}-\bm{H}_h\|} & \multicolumn{2}{|c|}{\|u-u_h\|} &  \multicolumn{2}{|c}{\|\bm{q}-\bm{q}_h\|} & \multicolumn{2}{|c|}{\|\bm{H}-\bm{H}_h\|}  \\
    p & n & \mbox{error}  & \mbox{order} & \mbox{error}  & \mbox{order} & \mbox{error}  & \mbox{order} & \mbox{error}  & \mbox{order} & \mbox{error}  & \mbox{order} & \mbox{error}  & \mbox{order}  \\
  \hline    
& 4  &  1.23\mbox{e-}4  &  --  &  2.30\mbox{e-}3  &  --  &  9.93\mbox{e-}2  &  --  &  1.23\mbox{e-}4  &  --  &  2.3\mbox{e-}3  &  --  &  9.93\mbox{e-}2  &  --  \\  
 & 8  &  1.96\mbox{e-}5  &  1.96  &  6.25\mbox{e-}4  &  1.88  &  5.05\mbox{e-}2  &   0.98  &  1.96\mbox{e-}5  &  1.96  &  6.25\mbox{e-}4  &  1.88  &  5.05\mbox{e-}2  &   0.98  \\  
1 & 16  &  1.99\mbox{e-}5  &  1.89  &  1.92\mbox{e-}4  &  1.70  &  2.54\mbox{e-}2  &   0.99  &  1.99\mbox{e-}5  &  1.89  &  1.92\mbox{e-}4  &  1.70  &  2.54\mbox{e-}2  &   0.99  \\  
 & 32  &  1.43\mbox{e-}5  &  1.54  &  8.47\mbox{e-}5  &  1.18  &  1.28\mbox{e-}2  &   0.99  &  1.43\mbox{e-}5  &  1.54  &  8.47\mbox{e-}5  &  1.18  &  1.28\mbox{e-}2  &   0.99  \\  
 & 64  &  8.35\mbox{e-}6  &  1.12  &  4.43\mbox{e-}5  &   0.94  &  6.40\mbox{e-}3  &  1.00  &  8.35\mbox{e-}6  &  1.12  &  4.43\mbox{e-}5  &   0.94  &  6.40\mbox{e-}3  &  1.00  \\  
\hline
 & 4  &  8.31\mbox{e-}6  &  --  &  9.82\mbox{e-}5  &  --  &  8.54\mbox{e-}3  &  --  &  8.31\mbox{e-}6  &  --  &  9.82\mbox{e-}5  &  --  &  8.54\mbox{e-}3  &  --  \\  
 & 8  &  2.57\mbox{e-}6  &  2.83  &  1.91\mbox{e-}5  &  2.36  &  2.24\mbox{e-}3  &  1.93  &  2.57\mbox{e-}6  &  2.83  &  1.91\mbox{e-}5  &  2.36  &  2.24\mbox{e-}3  &  1.93  \\  
2 & 16  &  7.00\mbox{e-}7  &  2.71  &  4.43\mbox{e-}6  &  2.11  &  5.69\mbox{e-}4  &  1.97  &  7.00\mbox{e-}7  &  2.71  &  4.43\mbox{e-}6  &  2.11  &  5.69\mbox{e-}4  &  1.97  \\  
 & 32  &  1.81\mbox{e-}7  &  2.42  &  1.10\mbox{e-}6  &  2.01  &  1.43\mbox{e-}4  &  1.99  &  1.81\mbox{e-}7  &  2.42  &  1.10\mbox{e-}6  &  2.01  &  1.43\mbox{e-}4  &  1.99  \\  
 & 64  &  4.57\mbox{e-}8  &  2.16  &  2.75\mbox{e-}7  &  2.00  &  3.60\mbox{e-}5  &  2.00  &  4.57\mbox{e-}8  &  2.16  &  2.75\mbox{e-}7  &  2.00  &  3.60\mbox{e-}5  &  2.00  \\  
\hline
 & 4  &  3.16\mbox{e-}7  &  --  &  5.56\mbox{e-}6  &  --  &  9.39\mbox{e-}4  &  --  &  3.16\mbox{e-}7  &  --  &  5.56\mbox{e-}6  &  --  &  9.39\mbox{e-}4  &  --  \\  
 & 8  &  1.01\mbox{e-}8  &  3.83  &  3.95\mbox{e-}7  &  3.81  &  1.28\mbox{e-}4  &  2.87  &  1.01\mbox{e-}8  &  3.83  &  3.95\mbox{e-}7  &  3.81  &  1.28\mbox{e-}4  &  2.87  \\  
3 & 16  &  1.57\mbox{e-}9  &  3.93  &  2.72\mbox{e-}8  &  3.86  &  1.65\mbox{e-}5  &  2.96  &  1.57\mbox{e-}9  &  3.93  &  2.72\mbox{e-}8  &  3.86  &  1.65\mbox{e-}5  &  2.96  \\  
 & 32  &  3.19\mbox{e-}10  &  3.91  &  2.46\mbox{e-}9  &  3.47  &  2.07\mbox{e-}6  &  2.99  &  3.2\mbox{e-}10  &  3.91  &  2.47\mbox{e-}9  &  3.46  &  2.07\mbox{e-}6  &  2.99  \\  
 & 64  &  5.28\mbox{e-}11  &  3.67  &  3.24\mbox{e-}10  &  2.92  &  2.60\mbox{e-}7  &  2.99  &  5.32\mbox{e-}11  &  3.67  &  3.33\mbox{e-}10  &  2.89  &  2.60\mbox{e-}7  &  2.99  \\  
\hline
     \end{array} $
}
\caption{Error and convergence rate for the HDG methods for the Example 2 with $R = 2$.}
\label{ex2tab1}
\end{center}{$\phantom{|}$}     
\end{table}

\begin{table}[thbp]
  \begin{center}
\scalebox{0.8}{%
    $\begin{array}{|c|c||c c  | c c|  c c ||  c c| c c| c c|}
    \hline    
    & & \multicolumn{6}{|c||}{\mbox{Newton-HDG method}} & \multicolumn{6}{|c|}{\mbox{Fixed point-HDG method}} \\ 
   \hline 
  \mbox{degree} & \mbox{mesh} & \multicolumn{2}{|c}{\|u-u_h\|} &  \multicolumn{2}{|c}{\|\bm{q}-\bm{q}_h\|} & \multicolumn{2}{|c||}{\|\bm{H}-\bm{H}_h\|} & \multicolumn{2}{|c|}{\|u-u_h\|} &  \multicolumn{2}{|c}{\|\bm{q}-\bm{q}_h\|} & \multicolumn{2}{|c|}{\|\bm{H}-\bm{H}_h\|}  \\
    p & n & \mbox{error}  & \mbox{order} & \mbox{error}  & \mbox{order} & \mbox{error}  & \mbox{order} & \mbox{error}  & \mbox{order} & \mbox{error}  & \mbox{order} & \mbox{error}  & \mbox{order}  \\
  \hline    
& 4  &  1.79\mbox{e-}3  &  --  &  1.14\mbox{e-}2  &  --  &  5.27\mbox{e-}1  &  --  &  1.79\mbox{e-}3  &  --  &  1.14\mbox{e-}2  &  --  &  5.27\mbox{e-}1  &  --  \\  
 & 8  &  4.49\mbox{e-}4  &  1.84  &  4.01\mbox{e-}3  &  1.50  &  3.46\mbox{e-}1  &   0.61  &  4.49\mbox{e-}4  &  1.84  &  4.01\mbox{e-}3  &  1.50  &  3.46\mbox{e-}1  &   0.61  \\  
1 & 16  &  6.87\mbox{e-}5  &  1.98  &  1.13\mbox{e-}3  &  1.83  &  2.01\mbox{e-}1  &   0.78  &  6.87\mbox{e-}5  &  1.98  &  1.13\mbox{e-}3  &  1.83  &  2.01\mbox{e-}1  &   0.78  \\  
 & 32  &  3.31\mbox{e-}5  &  1.82  &  3.13\mbox{e-}4  &  1.85  &  1.09\mbox{e-}1  &   0.89  &  3.31\mbox{e-}5  &  1.82  &  3.13\mbox{e-}4  &  1.85  &  1.09\mbox{e-}1  &   0.89  \\  
 & 64  &  2.70\mbox{e-}5  &  1.20  &  1.51\mbox{e-}4  &  1.05  &  5.64\mbox{e-}2  &   0.95  &  2.70\mbox{e-}5  &  1.20  &  1.51\mbox{e-}4  &  1.05  &  5.64\mbox{e-}2  &   0.95  \\  
\hline
 & 4  &  1.35\mbox{e-}4  &  --  &  1.76\mbox{e-}3  &  --  &  2.10\mbox{e-}1  &  --  &  1.35\mbox{e-}4  &  --  &  1.76\mbox{e-}3  &  --  &  2.10\mbox{e-}1  &  --  \\  
 & 8  &  6.99\mbox{e-}5  &  2.39  &  6.67\mbox{e-}4  &  1.40  &  8.49\mbox{e-}2  &  1.31  &  6.99\mbox{e-}5  &  2.39  &  6.67\mbox{e-}4  &  1.40  &  8.49\mbox{e-}2  &  1.31  \\  
2 & 16  &  2.04\mbox{e-}5  &  2.33  &  2.13\mbox{e-}4  &  1.65  &  2.77\mbox{e-}2  &  1.62  &  2.04\mbox{e-}5  &  2.33  &  2.13\mbox{e-}4  &  1.65  &  2.77\mbox{e-}2  &  1.62  \\  
 & 32  &  5.23\mbox{e-}6  &  2.18  &  5.75\mbox{e-}5  &  1.89  &  7.88\mbox{e-}3  &  1.81  &  5.23\mbox{e-}6  &  2.18  &  5.75\mbox{e-}5  &  1.89  &  7.88\mbox{e-}3  &  1.81  \\  
 & 64  &  1.32\mbox{e-}6  &  2.06  &  1.47\mbox{e-}5  &  1.96  &  2.09\mbox{e-}3  &  1.91  &  1.32\mbox{e-}6  &  2.06  &  1.47\mbox{e-}5  &  1.96  &  2.09\mbox{e-}3  &  1.91  \\  
\hline
 & 4  &  8.68\mbox{e-}5  &  --  &  6.98\mbox{e-}4  &  --  &  8.27\mbox{e-}2  &  --  &  8.68\mbox{e-}5  &  --  &  6.98\mbox{e-}4  &  --  &  8.27\mbox{e-}2  &  --  \\  
 & 8  &  8.09\mbox{e-}6  &  3.16  &  1.12\mbox{e-}4  &  2.64  &  2.28\mbox{e-}2  &  1.86  &  8.09\mbox{e-}6  &  3.16  &  1.12\mbox{e-}4  &  2.64  &  2.28\mbox{e-}2  &  1.86  \\  
3 & 16  &  5.52\mbox{e-}7  &  3.57  &  1.23\mbox{e-}5  &  3.18  &  4.37\mbox{e-}3  &  2.38  &  5.52\mbox{e-}7  &  3.57  &  1.23\mbox{e-}5  &  3.18  &  4.37\mbox{e-}3  &  2.38  \\  
 & 32  &  3.02\mbox{e-}8  &  3.83  &  9.56\mbox{e-}7  &  3.69  &  6.53\mbox{e-}4  &  2.74  &  3.02\mbox{e-}8  &  3.83  &  9.56\mbox{e-}7  &  3.69  &  6.53\mbox{e-}4  &  2.74  \\  
 & 64  &  1.33\mbox{e-}9  &  3.95  &  5.98\mbox{e-}8  &  4.00  &  8.66\mbox{e-}5  &  2.91  &  1.33\mbox{e-}9  &  3.95  &  5.98\mbox{e-}8  &  4.00  &  8.66\mbox{e-}5  &  2.91  \\  
\hline
     \end{array} $
}
\caption{Error and convergence rate for the HDG methods for the Example 2 with $R = \sqrt{2} + 0.1$.}
\label{ex2tab2}
\end{center}{$\phantom{|}$}     
\end{table}

\begin{table}[ht]
  \begin{center}
\scalebox{0.8}{%
    $\begin{array}{|c|c||c c  | c c|  c c ||  c c| c c| c c|}
    \hline    
    & & \multicolumn{6}{|c||}{\mbox{Newton-HDG method}} & \multicolumn{6}{|c|}{\mbox{Fixed point-HDG method}} \\ 
   \hline 
  \mbox{degree} & \mbox{mesh} & \multicolumn{2}{|c}{\|u-u_h\|} &  \multicolumn{2}{|c}{\|\bm{q}-\bm{q}_h\|} & \multicolumn{2}{|c||}{\|\bm{H}-\bm{H}_h\|} & \multicolumn{2}{|c|}{\|u-u_h\|} &  \multicolumn{2}{|c}{\|\bm{q}-\bm{q}_h\|} & \multicolumn{2}{|c|}{\|\bm{H}-\bm{H}_h\|}  \\
    p & n & \mbox{error}  & \mbox{order} & \mbox{error}  & \mbox{order} & \mbox{error}  & \mbox{order} & \mbox{error}  & \mbox{order} & \mbox{error}  & \mbox{order} & \mbox{error}  & \mbox{order}  \\
  \hline    
 & 16  &  9.17\mbox{e-}4  &  --  &  1.62\mbox{e-}2  &  --  &  1.8\mbox{e-}0  &  --  &  9.17\mbox{e-}4  &  --  &  1.62\mbox{e-}2  &  --  &  1.8\mbox{e-}0  &  --  \\  
 & 32  &  2.41\mbox{e-}4  &  1.85  &  7.71\mbox{e-}3  &  1.07  &  1.61\mbox{e-}0  &   0.16  &  2.41\mbox{e-}4  &  1.85  &  7.71\mbox{e-}3  &  1.07  &  1.61\mbox{e-}0  &   0.16  \\  
1 & 64  &  6.28\mbox{e-}5  &  1.88  &  3.01\mbox{e-}3  &  1.36  &  1.18\mbox{e-}0  &   0.44  &  6.28\mbox{e-}5  &  1.88  &  3.01\mbox{e-}3  &  1.36  &  1.18\mbox{e-}0  &   0.44  \\  
 & 128  &  2.45\mbox{e-}5  &  1.51  &  9.79\mbox{e-}4  &  1.62  &  7.35\mbox{e-}1  &   0.68  &  2.45\mbox{e-}5  &  1.51  &  9.79\mbox{e-}4  &  1.62  &  7.35\mbox{e-}1  &   0.68  \\  
 & 256  &  1.31\mbox{e-}5  &  1.04  &  2.83\mbox{e-}4  &  1.79  &  4.10\mbox{e-}1  &   0.84  &  1.31\mbox{e-}5  &  1.04  &  2.83\mbox{e-}4  &  1.79  &  4.10\mbox{e-}1  &   0.84  \\  
\hline
 & 16  &  9.57\mbox{e-}5  &  --  &  3.19\mbox{e-}3  &  --  &  1.21\mbox{e-}0  &  --  &  9.57\mbox{e-}5  &  --  &  3.19\mbox{e-}3  &  --  &  1.21\mbox{e-}0  &  --  \\  
 & 32  &  3.60\mbox{e-}5  &  2.45  &  1.33\mbox{e-}3  &  1.27  &  7.23\mbox{e-}1  &   0.74  &  3.60\mbox{e-}5  &  2.45  &  1.33\mbox{e-}3  &  1.27  &  7.23\mbox{e-}1  &   0.74  \\  
2 & 64  &  9.89\mbox{e-}6  &  2.67  &  5.37\mbox{e-}4  &  1.30  &  3.35\mbox{e-}1  &  1.11  &  9.89\mbox{e-}6  &  2.67  &  5.37\mbox{e-}4  &  1.30  &  3.35\mbox{e-}1  &  1.11  \\  
 & 128  &  2.44\mbox{e-}6  &  2.43  &  1.56\mbox{e-}4  &  1.78  &  1.21\mbox{e-}1  &  1.47  &  2.44\mbox{e-}6  &  2.43  &  1.56\mbox{e-}4  &  1.78  &  1.21\mbox{e-}1  &  1.47  \\  
 & 256  &  6.03\mbox{e-}7  &  2.15  &  4.02\mbox{e-}5  &  1.96  &  3.66\mbox{e-}2  &  1.73  &  6.03\mbox{e-}7  &  2.15  &  4.02\mbox{e-}5  &  1.96  &  3.66\mbox{e-}2  &  1.73  \\  
\hline
 & 16  &  1.04\mbox{e-}4  &  --  &  3.95\mbox{e-}3  &  --  &  6.84\mbox{e-}1  &  --  &  1.04\mbox{e-}4  &  --  &  3.95\mbox{e-}3  &  --  &  6.84\mbox{e-}1  &  --  \\  
 & 32  &  1.05\mbox{e-}5  &  3.12  &  9.06\mbox{e-}4  &  2.12  &  3.28\mbox{e-}1  &  1.06  &  1.05\mbox{e-}5  &  3.12  &  9.06\mbox{e-}4  &  2.12  &  3.28\mbox{e-}1  &  1.06  \\  
3 & 64  &  8.19\mbox{e-}7  &  3.31  &  1.39\mbox{e-}4  &  2.70  &  1.06\mbox{e-}1  &  1.63  &  8.19\mbox{e-}7  &  3.31  &  1.39\mbox{e-}4  &  2.70  &  1.06\mbox{e-}1  &  1.63  \\  
 & 128  &  5.29\mbox{e-}8  &  3.53  &  1.40\mbox{e-}5  &  3.31  &  2.29\mbox{e-}2  &  2.21  &  5.29\mbox{e-}8  &  3.53  &  1.40\mbox{e-}5  &  3.31  &  2.28\mbox{e-}2  &  2.21  \\  
 & 256  &  3.22\mbox{e-}9  &  3.77  &  1.00\mbox{e-}6  &  3.80  &  3.63\mbox{e-}3  &  2.65  &  3.22\mbox{e-}9  &  3.77  &  1.00\mbox{e-}6  &  3.80  &  3.63\mbox{e-}3  &  2.65  \\  
\hline
     \end{array} $
}
\caption{Error and convergence rate for the HDG methods for the Example 2 with $R = \sqrt{2} + 0.01$.}
\label{ex2tab3}
\end{center}{$\phantom{|}$}     
\end{table}

\begin{table}[ht]
  \begin{center}
\scalebox{0.8}{%
    $\begin{array}{|c||c c c | c c c|  c c c ||  c c c| c c c| c c c|}
    \hline    
    & \multicolumn{9}{|c||}{\mbox{Newton-HDG method}} & \multicolumn{9}{|c|}{\mbox{Fixed point-HDG method}} \\ 
   \hline 
   \mbox{mesh} & \multicolumn{3}{|c}{R=2} &  \multicolumn{3}{|c}{R=\sqrt{2}+0.1} & \multicolumn{3}{|c||}{R=\sqrt{2}+0.01} & \multicolumn{3}{|c|}{R=2} &  \multicolumn{3}{|c}{R=\sqrt{2}+0.1} & \multicolumn{3}{|c|}{R=\sqrt{2}+0.01}  \\
     \mbox{level} & \mbox{P1}  & \mbox{P2} & \mbox{P3}  & \mbox{P1}  & \mbox{P2} & \mbox{P3} & \mbox{P1}  & \mbox{P2} & \mbox{P3}  & \mbox{P1}  & \mbox{P2} & \mbox{P3}  & \mbox{P1}  & \mbox{P2} & \mbox{P3} & \mbox{P1}  & \mbox{P2} & \mbox{P3}   \\
  \hline    
 1  &  6  &  6  &  6  &  6  &  6  &  6  &  7  &  7  &  7  &  27  &  30  &  31  &  29  &  34  &  39  &  35  &  59  &  72  \\  
 2  &  6  &  6  &  6  &  6  &  6  &  6  &  7  &  7  &  7  &  29  &  31  &  32  &  33  &  38  &  43  &  39  &  54  &  84  \\  
 3  &  6  &  6  &  6  &  7  &  6  &  6  &  7  &  7  &  7  &  30  &  32  &  33  &  37  &  42  &  43  &  47  &  66  &  97  \\  
 4  &  7  &  6  &  6  &  7  &  6  &  6  &  7  &  7  &  7  &  32  &  33  &  34  &  40  &  44  &  46  &  54  &  89  &  121  \\  
 5  &  7  &  6  &  6  &  7  &  6  &  6  &  7  &  7  &  7  &  33  &  34  &  35  &  43  &  46  &  49  &  60  &  112  &  146  \\  
\hline
     \end{array} $
}
\caption{Number of iterations required to reach convergence for the HDG methods for the Example 2.}
\label{ex2tab4}
\end{center}{$\phantom{|}$}     
\end{table}

\section{Optimal transport for $r$-adaptive mesh generation}

In this section, we  review optimal transport theory  and describe how it can be applied to $r$-adaptive mesh generation. This $r$-adaptive mesh generation approach results in the Monge–Amp\`ere equation with nonlinear Neumann boundary condition, which is solved by extending the HDG methods described in the previous section. We present numerical experiments to demonstrate the performance of the HDG methods for $r$-adaptive mesh generation.

\subsection{Optimal transport theory}

The optimal transport (OT) problem is described as follows. Suppose we are given two probability densities: $\rho(\bm x)$ supported on $\Omega \in \mathbb{R}^d$ and $\rho'(\bm x')$ supported on $\Omega' \in \mathbb{R}^d$. The source density $\rho(\bm x)$ may be discontinuous and  even vanish. The target density $\rho'(\bm x')$ must be strictly positive and Lipschitz continuous. The OT problem is to find a map $\bm \phi : \Omega \to \Omega'$ such that it minimizes the following functional
\begin{equation}
\inf_{\bm \phi \in \mathcal{M}} \int_{\Omega} \|\bm x - \bm \phi(\bm x)\|^2 \rho(\bm x) d \bm x ,     
\end{equation}
where 
\begin{equation}
\label{otset}
\mathcal{M} = \{\bm \phi : \Omega \to \Omega', \ \rho'(\bm \phi(\bm x)) \det (\nabla \bm \phi(\bm x)) = \rho(\bm x), \ \forall \bm x \in \Omega  \} ,
\end{equation}
is the set of mappings which map the source density $\rho(\bm x)$ onto the target density $\rho'(\bm x')$.  

In \cite{Brenier1991}, Brenier gave the proof of the existence and uniqueness of the solution of the OT problem. Furthermore, the optimal map $\bm \phi$ can be written as the gradient  of a unique (up to a constant) convex potential $u$, so that $\bm \phi(\bm x) = \nabla u(\bm x)$, $\Delta u(\bm x) > 0$. Substituting $\bm \phi(\bm x) = \nabla u(\bm x)$ into (\ref{otset}) results in the Monge–Amp\`ere equation
\begin{equation}
\label{mae}
 \rho'(\nabla u(\bm x)) \det (D^2 u(\bm x)) = \rho(\bm x) \quad \mbox{in } \Omega, 
\end{equation}
along with the restriction that $u$ is convex.  The equation lacks standard boundary conditions. However, it is geometrically constrained by the fact that the gradient map takes $\partial \Omega$ to $\partial \Omega'$: 
\begin{equation}
\label{maebc}
 \nabla u(\bm x) \in \partial \Omega', \quad \forall \bm x \in \Omega .
\end{equation}
This constraint is referred to as the second boundary value problem for the Monge–Amp\`ere equation. In \cite{Benamou2014}, the constraint (\ref{maebc}) is replaced with a Hamilton–Jacobi equation on the boundary
\begin{equation}
\label{maehj}
 H(\nabla u(\bm x)) := \mbox{dist} (\nabla u(\bm x), \partial \Omega') = 0, \quad \forall \bm x \in \partial \Omega .
\end{equation}
If the boundary $\partial \Omega'$ can be expressed by $\partial \Omega' = \{\bm x' \in \Omega' : g(\bm x') = 0\}$ then the  Hamilton–Jacobi equation reduces to the following Neumann boundary condition
\begin{equation}
\label{maenm}
 g(\nabla u(\bm x)) = 0, \quad \forall \bm x \in \partial \Omega .
\end{equation}
This boundary condition simplifies to a linear Neumann boundary condition when $g(\cdot)$ is a linear function, that is, when the boundary $\partial \Omega'$ is flat. For certain problems where densities are periodic, it is natural and convenient to use periodic boundary conditions instead. 

\subsection{Adaptive mesh redistribution}

One approach to mesh adaptation is the equidistribution principle that equidistributes the target density function $\rho'$ so that the source density $\rho$ is uniform on $\Omega$ \cite{Delzanno2008,Chacon2011}. The equidistribution principle leads to a constant source density $\rho(\bm x) = \theta$, where the constant $\theta$ is given by $\theta = \int_{\Omega'} \rho'(\bm x') d \bm x' / \int_{\Omega} d \bm x$. Using the optimal transport theory, the optimal mesh is sought by solving the  Monge–Amp\`ere equation with the Neumann boundary condition:
\begin{equation}
\label{maem}
\begin{split}
  \det (D^2 u(\bm x)) & = f(\nabla u(\bm x)), \quad \mbox{in } \Omega, \\
 g(\nabla u(\bm x)) & = 0, \qquad \qquad \mbox{on } \partial \Omega , \\
 \int_{\Omega} u(\bm x) d \bm x &= 0,
 \end{split}
\end{equation}
where $f(\nabla u(\bm x)) = \theta/\rho'(\nabla u(\bm x))$. The gradient of $u$ gives us the desired mesh.


In two dimensions, we can rewrite the above equation as a first-order system of equations
\begin{equation}
\begin{array}{rcll}
\bm{H} - \nabla \bm q & = & 0, \quad & \mbox{in } \Omega, \\
\bm{q} - \nabla u & = & 0, \quad & \mbox{in } \Omega, \\
s(\bm H, \bm q) - \nabla \cdot \bm q & = & 0 , \quad & \mbox{in }  \Omega, \\
g(\bm q) & = & 0, \quad & \mbox{on } \partial \Omega , \\
\int_{\Omega} u(\bm x) d \bm x &= & 0,
\end{array}
\label{maem3}
\end{equation}
where $s(\bm H,\bm q) = \sqrt{H_{11}^2 + H_{22}^2 + H_{12}^2 + H_{21}^2 + 2 f(\bm q)}$. The system (\ref{maem3}) differs from (\ref{model_problem3}) in the source term and the boundary condition. In order to solve (\ref{maem3}), we will need to extend the HDG methods described in the previous section.

\subsection{The Newton-HDG formulation}

The HDG discretization of the system (\ref{maem3}) is to find $(\bm{H}_h, \bm{q}_h,u_h,\widehat{u}_h) \in \bm{W}_{h}^p \times \bm{V}_{h}^p  \times U_{h}^p \times M_h^p$ such that
\begin{equation}
\label{eq82}
\begin{array}{rcl}
 \left(\bm{H}_h, \bm{G}\right)_{\mathcal{T}_h} +  \left(\bm q_h,  \nabla \cdot \bm{G}\right)_{\mathcal{T}_h} -  \left\langle \widehat{\bm q}_h, \bm{G} \cdot \bm{n} \right\rangle_{\partial \mathcal{T}_h} & = & 0,  \\
 \left(\bm{q}_h, \bm{v}\right)_{\mathcal{T}_h} +  \left(u_h,  \nabla \cdot \bm{v}\right)_{\mathcal{T}_h} -  \left\langle \widehat{u}_h, \bm{v} \cdot \bm{n} \right\rangle_{\partial \mathcal{T}_h} & = & 0,  \\
\left(\bm{q}_h, \nabla w\right)_{\mathcal{T}_h}  -  \left\langle \widehat{\bm{q}}_h \cdot \bm{n}, w \right\rangle_{\partial \mathcal{T}_h} + (s(\bm H_h, \bm q_h),w)_{\mathcal{T}_h} & = & 0, \\
\left\langle \widehat{\bm{q}}_h  \cdot \bm{n} , \mu \right\rangle_{\partial \mathcal{T}_h \backslash \partial \Omega} + \left\langle g(\bm q_h) + \tau (\widehat{u}_h - u_h), \mu \right\rangle_{\partial \Omega} & =  & 0, \\
(u_h,1)_{\mathcal{T}_h} & = & 0,
\end{array}
\end{equation}
for all $(\bm G, \bm{v}, w, \mu) \in  \bm{W}_h^p \times 
 \bm{V}_h^p \times U_h^p \times M_h^p$, where
\begin{equation}
\widehat{\bm{q}}_h = {\bm{q}_h} -
\tau (u_h - \widehat{u}_h) \bm{n}, \quad \mbox{on } \mathcal{E}_h.
\label{fluxdef2}
\end{equation}
We are going to use Newton's method to solve this nonlinear system of equations.

For each Newton step $l$, we compute $\left(\delta \bm{H}_h^l, \delta  \bm{q}_h^l, \delta u_h^l, \delta 
 \widehat{u}_h^l \right) \in \bm{W}_{h}^p \times \bm{V}_{h}^p  \times U_{h}^p \times M_h^p$ by solving the following linear system
\begin{subequations}\label{neq:HDG_increments}
\begin{alignat}{1}
\left(\delta \bm{H}_h^l, \bm{G}\right)_{\mathcal{T}_h} +  \left(\delta \bm q_h^l,  \nabla \cdot \bm{G}\right)_{\mathcal{T}_h} -  \left\langle \delta \widehat{\bm q}_h^l, \bm{G} \cdot \bm{n} \right\rangle_{\partial \mathcal{T}_h}  & =  r_1(\bm G) , \label{nsubeq:increm1} \\
  \left(\delta \bm{q}_h^l, \bm{v}\right)_{\mathcal{T}_h} +  \left(\delta u_h^l,  \nabla \cdot \bm{v}\right)_{\mathcal{T}_h} -  \left\langle \delta \widehat{u}_h^l, \bm{v} \cdot \bm{n} \right\rangle_{\partial \mathcal{T}_h}  & = r_2(\bm{v}) , \label{nsubeq:increm2} \\
  (\partial s_{\bm H}(\bm H_h^l, \bm q_h^l) \delta \bm H_h^l ,w)_{\mathcal{T}_h} + (\partial s_{\bm q}(\bm H_h^l, \bm q_h^l) \delta \bm q_h^l ,w)_{\mathcal{T}_h} \nonumber \\ 
  + \left(\delta \bm{q}_h^l, \nabla w\right)_{\mathcal{T}_h}  -  \left\langle \delta \widehat{\bm{q}}_h^l \cdot \bm{n}, w \right\rangle_{\partial \mathcal{T}_h}  & = r_3(w) , \label{nsubeq:increm3} \\
  \left\langle \delta \widehat{\bm{q}}_h^l  \cdot \bm{n} , \mu \right\rangle_{\partial \mathcal{T}_h \backslash \partial \Omega}  + \left\langle \partial g_{\bm q}(\bm q_h^l) \cdot \delta \bm q_h^l + \tau (\delta \widehat{u}_h^l - \delta u_h^l) , \mu \right\rangle_{\partial \Omega}  & = r_4({\mu}), \label{nsubeq:increm4} \\
 (\delta u_h^l,1)_{\mathcal{T}_h} & =  r_5
\end{alignat}
\end{subequations}
for all $(\bm G, \bm{v}, w, \mu) \in  \bm{W}_h^p \times 
 \bm{V}_h^p \times U_h^p \times M_h^p$, the right-hand side residuals are given by
\begin{subequations}\label{neq:HDG_residuals}
\begin{alignat}{1}
  r_1(\bm G)     & = - \left(\bm{H}_h^l, \bm{G}\right)_{\mathcal{T}_h} -  \left(\bm q^l_h,  \nabla \cdot \bm{G}\right)_{\mathcal{T}_h} +  \left\langle \widehat{\bm q}_h^l, \bm{G} \cdot \bm{n} \right\rangle_{\partial \mathcal{T}_h} \\
  r_2(\bm{v})  & =- \left(\bm{q}_h^l, \bm{v}\right)_{\mathcal{T}_h} -  \left(u_h^l,  \nabla \cdot \bm{v}\right)_{\mathcal{T}_h} +  \left\langle \widehat{u}_h^l, \bm{v} \cdot \bm{n} \right\rangle_{\partial \mathcal{T}_h} \\
  r_3(w)& = - \left(\bm{q}_h^l, \nabla w\right)_{\mathcal{T}_h}  +  \left\langle \widehat{\bm{q}}_h^l \cdot \bm{n}, w \right\rangle_{\partial \mathcal{T}_h} - (s(\bm H_h^l, \bm q_h^l),w)_{\mathcal{T}_h}  \\
r_4(\mu) & = -\left\langle \widehat{\bm{q}}_h^l  \cdot \bm{n} , \mu \right\rangle_{\partial \mathcal{T}_h \backslash \partial \Omega} - \left\langle g(\bm q_h^l) + \tau(\widehat{u}_h^l - u_h^l), \mu \right\rangle_{\partial \Omega} \\
r_5 & = -( u_h^l,1)_{\mathcal{T}_h} .
\end{alignat}
\end{subequations}
Here $\partial s_{\bm H}(\bm H, \bm q)$ and $\partial s_{\bm q}(\bm H, \bm q)$ denote the partial derivatives of $s(\bm H, \bm q)$ with respect to $\bm H$ and $\bm q$, respectively. And $\partial g_{\bm q}(\bm q)$ denotes the partial derivatives of $g(\bm q)$ with respect to $\bm q$.

At each step of the Newton method, the linearization \eqref{neq:HDG_increments} gives the following matrix system to be solved
\begin{equation}\label{neq:matrix_system}
  \begin{pmatrix} \mathbb{A}^l & \mathbb{B}^l \\ \mathbb{C}^l & \mathbb{D}^l \\ \mathbb{E}^l & 0 \end{pmatrix} 
  \begin{pmatrix} \delta U^l  \\ \delta \widehat{U}^l \end{pmatrix}
= \begin{pmatrix} R_{123}^l  \\ R_4^l \\ r_5 \end{pmatrix} \\ 
\end{equation}
where $\delta U^l$ and $\delta \widehat{U}^l$ are the vectors of degrees of freedom of $(\delta \bm{H}_h^l, \delta  \bm{q}_h^l, \delta u_h^l)$ and $\delta 
 \widehat{u}_h^l$, respectively. The system \eqref{neq:matrix_system} is first solved for the traces only $\delta \widehat{U}^l$: 
\begin{equation}\label{neq:reduc_matrix_system}
\begin{pmatrix} \mathbb{K}^l \\ \mathbb{M}^l \end{pmatrix}
\delta \widehat{U}^l = \begin{pmatrix} R^l \\ r_5 \end{pmatrix}
\end{equation}
where 
\begin{equation}\label{neq:SchurCompl}
  \mathbb{K}^l = \mathbb{D}^l - \mathbb{C}^l \left( \mathbb{A}^l \right)^{-1} \mathbb{B}^l , \qquad \mathbb{M}^l =  - \mathbb{E}^l \left( \mathbb{A}^l \right)^{-1} \mathbb{B}^l, \qquad 
  R^l = R_4^l - \mathbb{C}^l \left( \mathbb{A}^l \right)^{-1} R_{123}^l .
\end{equation}
Once $\delta \widehat{U}^l$ is known, the other unknowns $\delta U^l$ are then retrieved element-wise.

\subsection{The fixed point-HDG formulation}

To devise the fixed point-HDG method for solving (\ref{maem3}), we describe how to deal with the boundary condition $g(\bm q) = 0$. For any given $\bm q^{\ell -1}$, we linearize it around $\bm q^{\ell -1}$ to obtain 
\begin{equation}
g(\bm q^{l-1}) + \partial g_{\bm q}(\bm q^{l-1}) \cdot \left( \bm q^{l} - \bm q^{l-1} \right) = 0 .
\end{equation}
Starting from an initial guess $(\bm H^0_h, \bm q_h^0, u_h^0)$  we  find $(\bm{q}_h^l,u_h^l,\widehat{u}_h^l) \in \bm{V}_{h}^p  \times U_{h}^p \times M_h^p$ such that
\begin{equation}
\label{nfpHDG}
\begin{array}{rcl}
  \left(\bm{q}_h^l, \bm{v}\right)_{\mathcal{T}_h} +  \left(u_h^l,  \nabla \cdot \bm{v}\right)_{\mathcal{T}_h} -  \left\langle \widehat{u}_h^l, \bm{v} \cdot \bm{n} \right\rangle_{\partial \mathcal{T}_h} & = & 0,  \\
\left(\bm{q}_h^l, \nabla w\right)_{\mathcal{T}_h}  -  \left\langle \widehat{\bm{q}}_h^l \cdot \bm{n}, w \right\rangle_{\partial \mathcal{T}_h}  & = & - (s(\bm H_h^{l-1}, \bm q_h^{l-1}),w)_{\mathcal{T}_h}, \\
\left\langle \widehat{\bm{q}}_h^l  \cdot \bm{n} , \mu \right\rangle_{\partial \mathcal{T}_h \backslash \partial \Omega}  + \left\langle \partial g_{\bm q}(\bm q_h^{l-1}) \cdot \bm q_h^l + \tau ( \widehat{u}_h^l -  u_h^l) , \mu \right\rangle_{\partial \Omega} & =  &  -\left\langle  b(\bm q^{l-1}) , \mu \right\rangle_{\partial \Omega}, \\
(u_h^l,1)_{\mathcal{T}_h} & = & 0, 
\end{array}
\end{equation}
for all $(\bm{v}, w, \mu) \in  \bm{V}_h^p \times U_h^p \times M_h^p$, and then compute $\bm{H}_h^l\in \bm{W}_{h}^p$ such that
\begin{equation}
\label{neqH}
 \left(\bm{H}_h^l, \bm{G}\right)_{\mathcal{T}_h} =  - \left(\bm q_h^l,  \nabla \cdot \bm{G}\right)_{\mathcal{T}_h} +  \left\langle \widehat{\bm q}_h^l, \bm{G} \cdot \bm{n} \right\rangle_{\partial \mathcal{T}_h}, \quad \forall \bm G \in \bm{W}_h^p .
\end{equation}
Note here that $b(\bm q_h^{l-1}) = g(\bm q_h^{l-1}) - \partial g_{\bm q}(\bm q_h^{l-1}) \cdot  \bm q_h^{l-1}$, and that the numerical flux $\widehat{\bm q}_h^l$ is defined by (\ref{fluxdef2}).

At each step of the fixed-point HDG method, the weak formulation (\ref{nfpHDG}) yields the matrix system similar to (\ref{neq:matrix_system}). The global matrix for the degrees of freedom of $\widehat{u}_h^l$ is changed at each step because of the linearization of the nonlinear boundary condition $g(\bm q) = 0$. In this case, the fixed-point HDG method is no longer competitive to the Newton-HDG method. We note however that if $g(\cdot)$ is a linear function, then the global matrix will be unchanged at each fixed-point step. In this case, the global matrix can be formed prior to performing the fixed-point iteration. In any case, the Newton-HDG method is more efficient than the fixed-point HDG method since the former requires considerably less number of iterations than the latter. 

\subsection{Numerical experiments}

We give several examples of high-order meshes generated using the HDG methods  for analytical density functions on both planar and curved domains. We compare the convergence of the HDG methods for these examples. We also demonstrate that our methods can generate smooth high-order meshes even at very high mesh resolutions.

\

\begin{figure}[hthbp]
	\centering
	\begin{subfigure}[b]{0.32\textwidth}
		\centering
		\includegraphics[width=\textwidth]{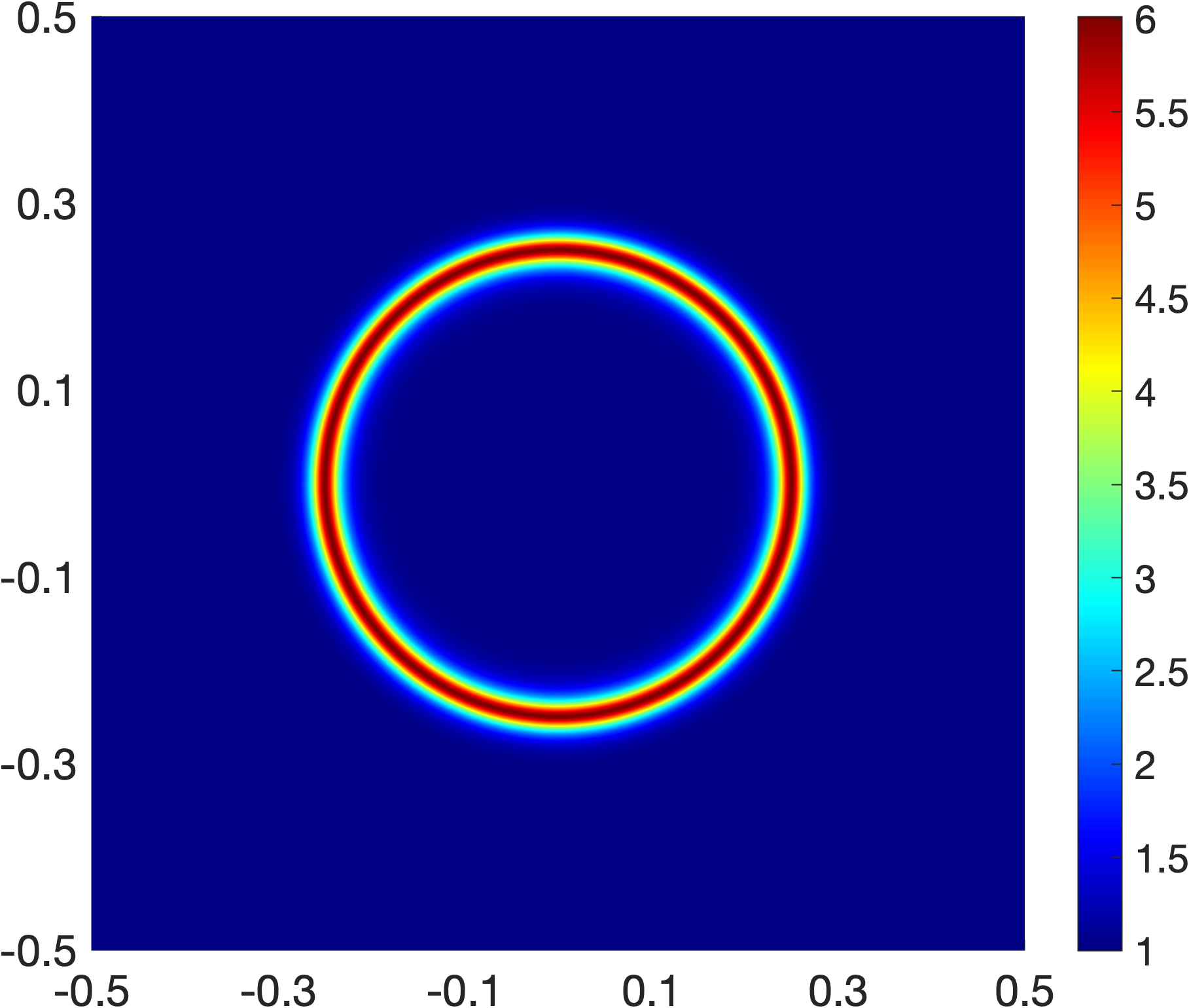}
		\caption{$\bm \alpha_1 = (5, 200, 0.25)$}
	\end{subfigure}
	\hfill
	\begin{subfigure}[b]{0.32\textwidth}
		\centering
		\includegraphics[width=\textwidth]{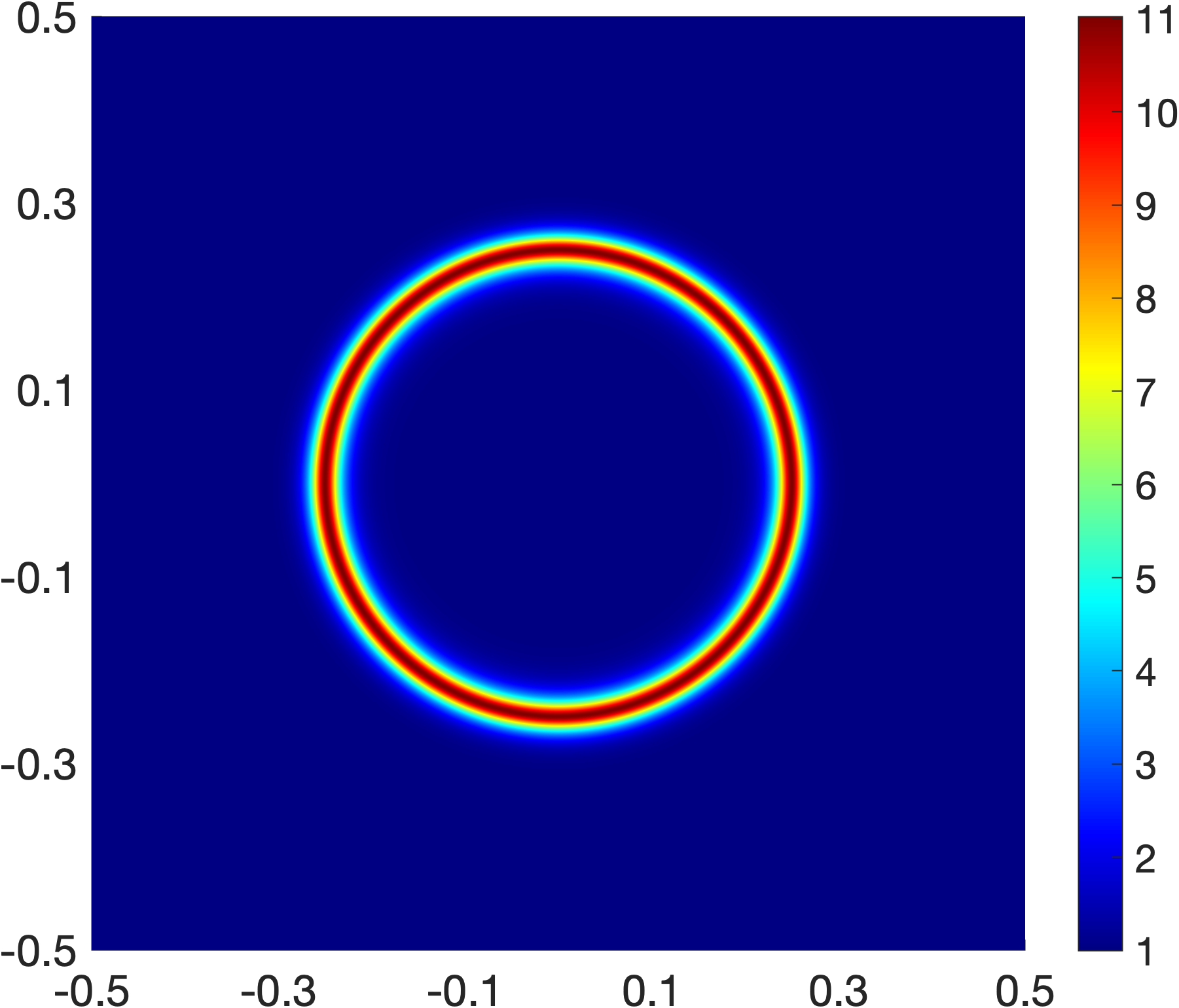}
		\caption{$\bm \alpha_2 = (10, 200, 0.25)$}
	\end{subfigure}
        \hfill
	\begin{subfigure}[b]{0.32\textwidth}
		\centering
		\includegraphics[width=\textwidth]{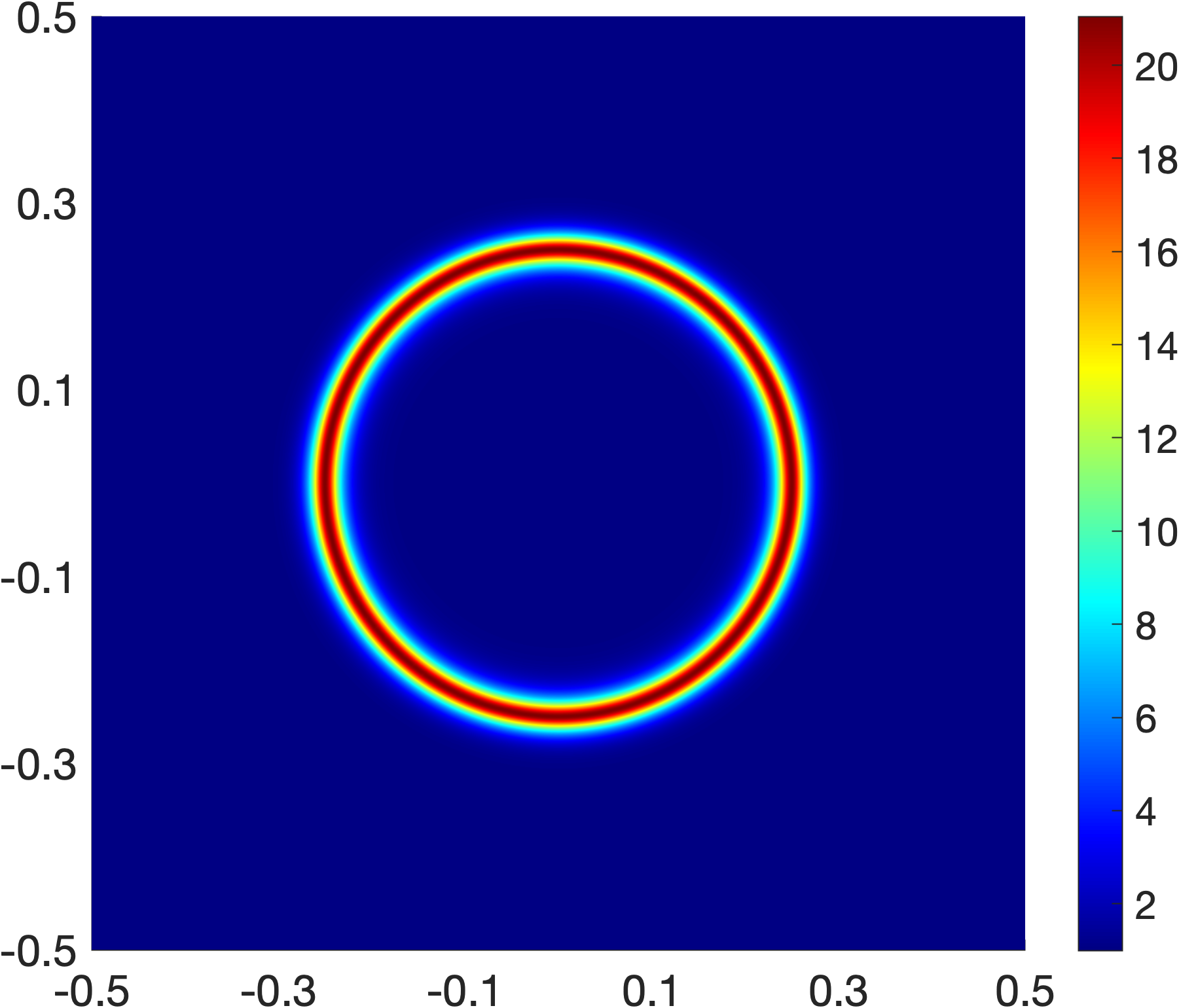}
		\caption{$\bm \alpha_3 = (20, 200, 0.25)$}
	\end{subfigure}
	\caption{Three instances of the density function (\ref{dens}).}
	\label{ex3fig1}
\end{figure}

\noindent 
\textbf{Ring meshes on a square domain.} We wish to  generate meshes on a square domain $\Omega' = (-0.5, 0.5)^2$ with the following target density function:
\begin{equation}
\label{dens}
\rho'(x', y') = 1 + \alpha_1 \mbox{sech}^2 \left( \alpha_2 \left( x'^2 + y'^2 - \alpha_3^2 \right) \right)  ,  
\end{equation}
where $\bm \alpha = (\alpha_1, \alpha_2, \alpha_3)$ determines the density function. This density function was introduced in \cite{Budd2015}. We consider three instances of the density function (\ref{dens}) corresponding to $\bm \alpha_1 = (5, 200, 0.25)$, $\bm \alpha_2 = (10, 200, 0.25)$, and $\bm \alpha_3 = (20, 200, 0.25)$, as shown in Figure \ref{ex3fig1}.

\begin{figure}[hthbp]
	\centering
	\begin{subfigure}[b]{0.32\textwidth}
		\centering
		\includegraphics[width=\textwidth]{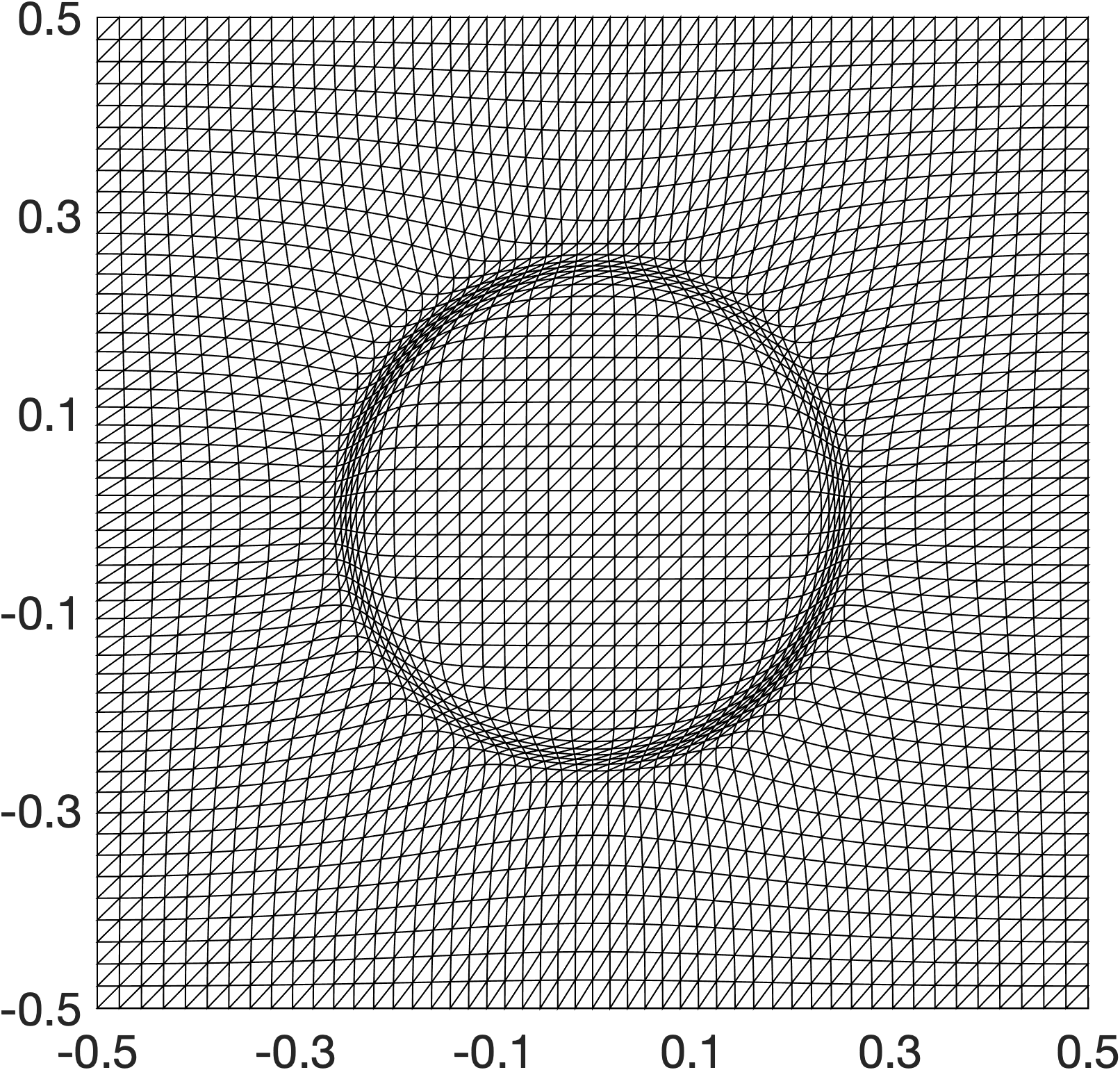}
		\caption{$\bm \alpha_1 = (5, 200, 0.25)$}
	\end{subfigure}
	\hfill
	\begin{subfigure}[b]{0.32\textwidth}
		\centering
		\includegraphics[width=\textwidth]{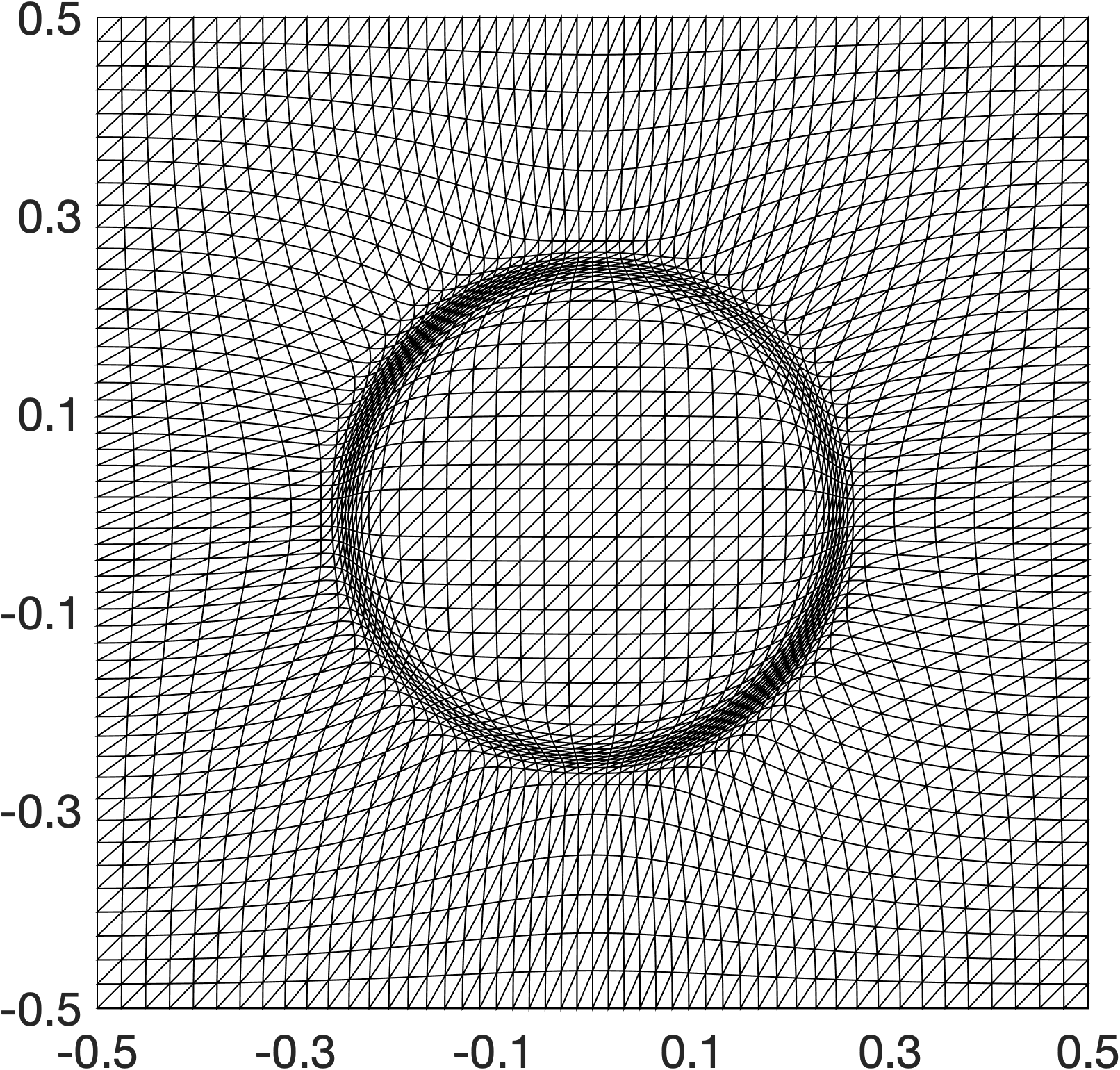}
		\caption{$\bm \alpha_2 = (10, 200, 0.25)$}
	\end{subfigure}
        \hfill
	\begin{subfigure}[b]{0.32\textwidth}
		\centering
		\includegraphics[width=\textwidth]{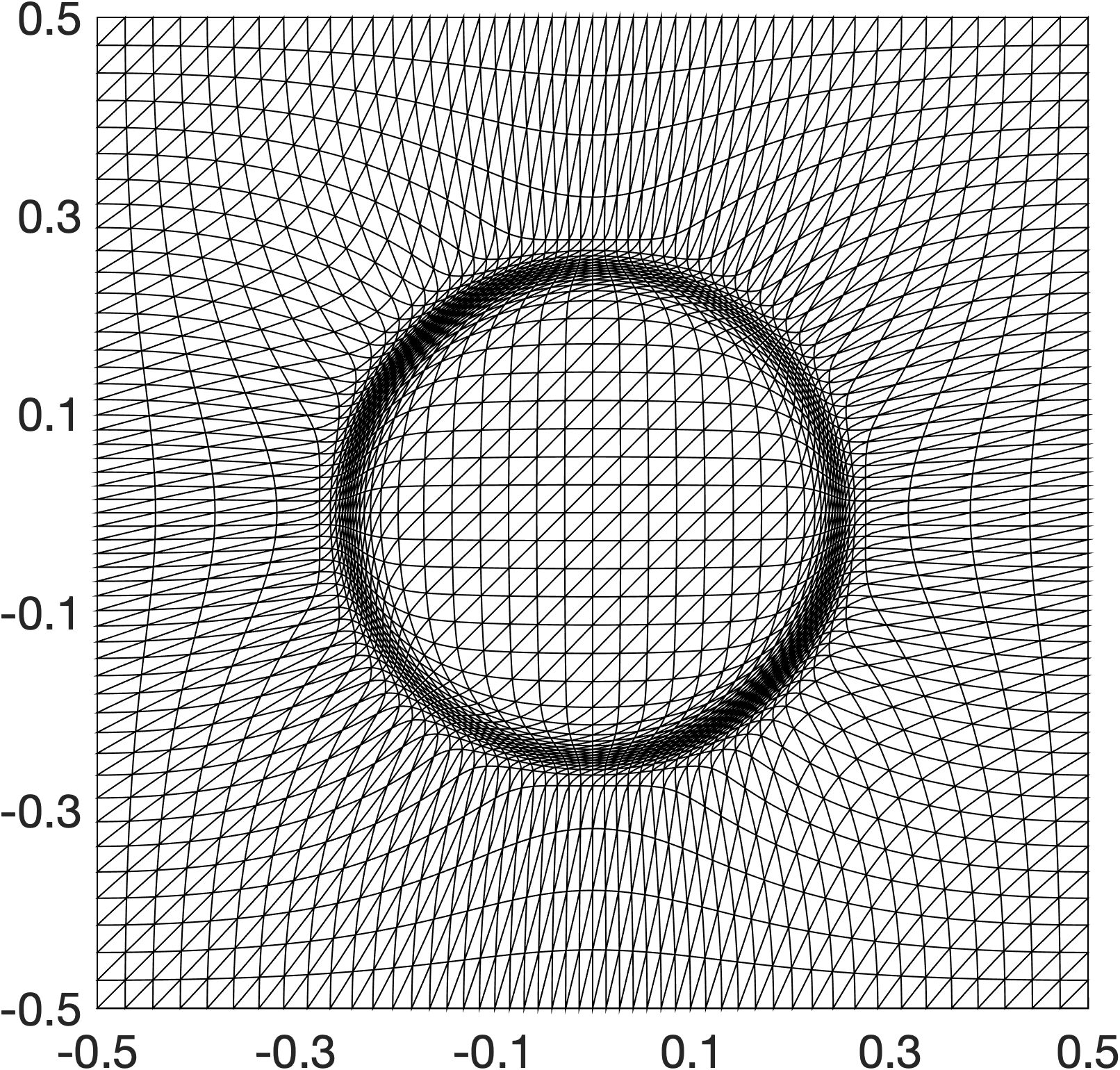}
		\caption{$\bm \alpha_3 = (20, 200, 0.25)$}
	\end{subfigure}
	\caption{Three high-order meshes are generated for the density functions shown in Figure \ref{ex3fig1} using a uniform background mesh of $50 \times 50 \times 2$ triangles with $p=3$.}
	\label{ex3fig2}
\end{figure}

First, a background mesh on which the optimal transport meshes are generated is a uniform grid of $50 \times 50 \times 2$ triangles on the square domain $\Omega = (-0.5, 0.5)^2$. Here we use polynomial degree $p=3$ to represent the numerical solution. Figure \ref{ex3fig2} depicts the three high-order meshes generated for these density instances. We observe that as $\alpha_1$ increases from 5, 10, to 20, it results in more elements concentrating into the ring.

\begin{figure}[hthbp]
	\centering
	\begin{subfigure}[b]{0.32\textwidth}
		\centering
		\includegraphics[width=\textwidth]{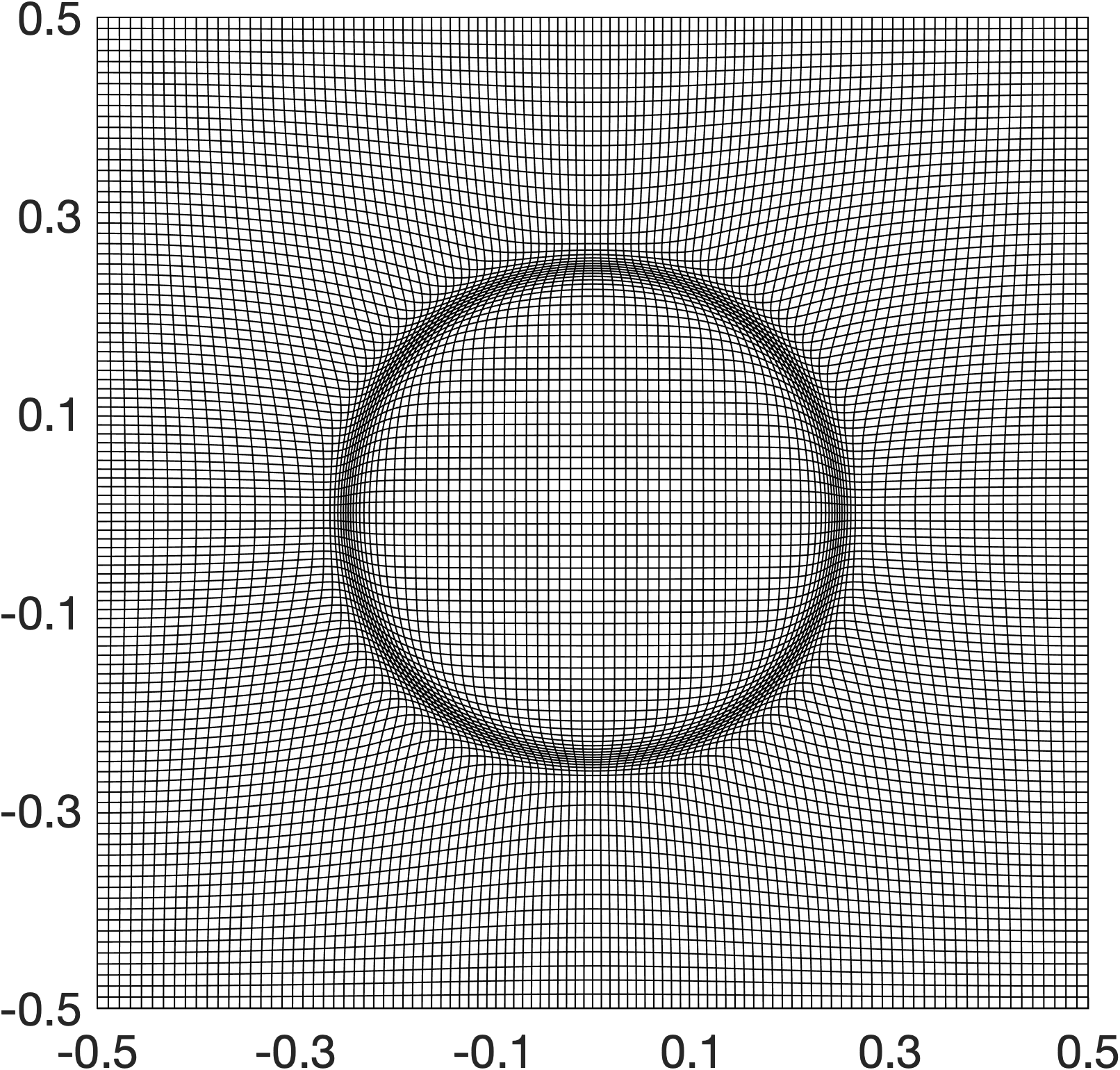}
		\caption{$\bm \alpha_1 = (5, 200, 0.25)$}
	\end{subfigure}
	\hfill
	\begin{subfigure}[b]{0.32\textwidth}
		\centering
		\includegraphics[width=\textwidth]{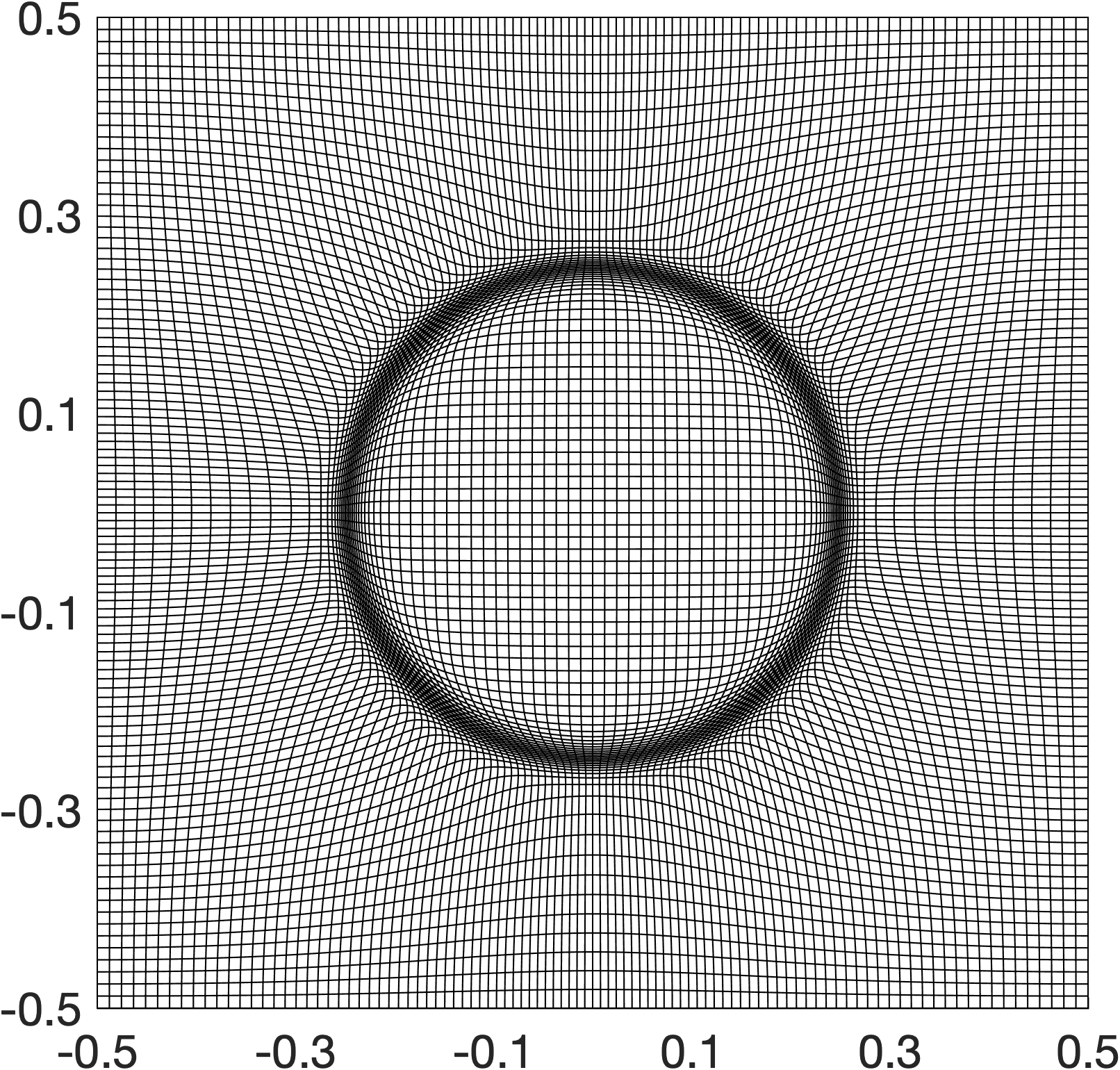}
		\caption{$\bm \alpha_2 = (10, 200, 0.25)$}
	\end{subfigure}
        \hfill
	\begin{subfigure}[b]{0.32\textwidth}
		\centering
		\includegraphics[width=\textwidth]{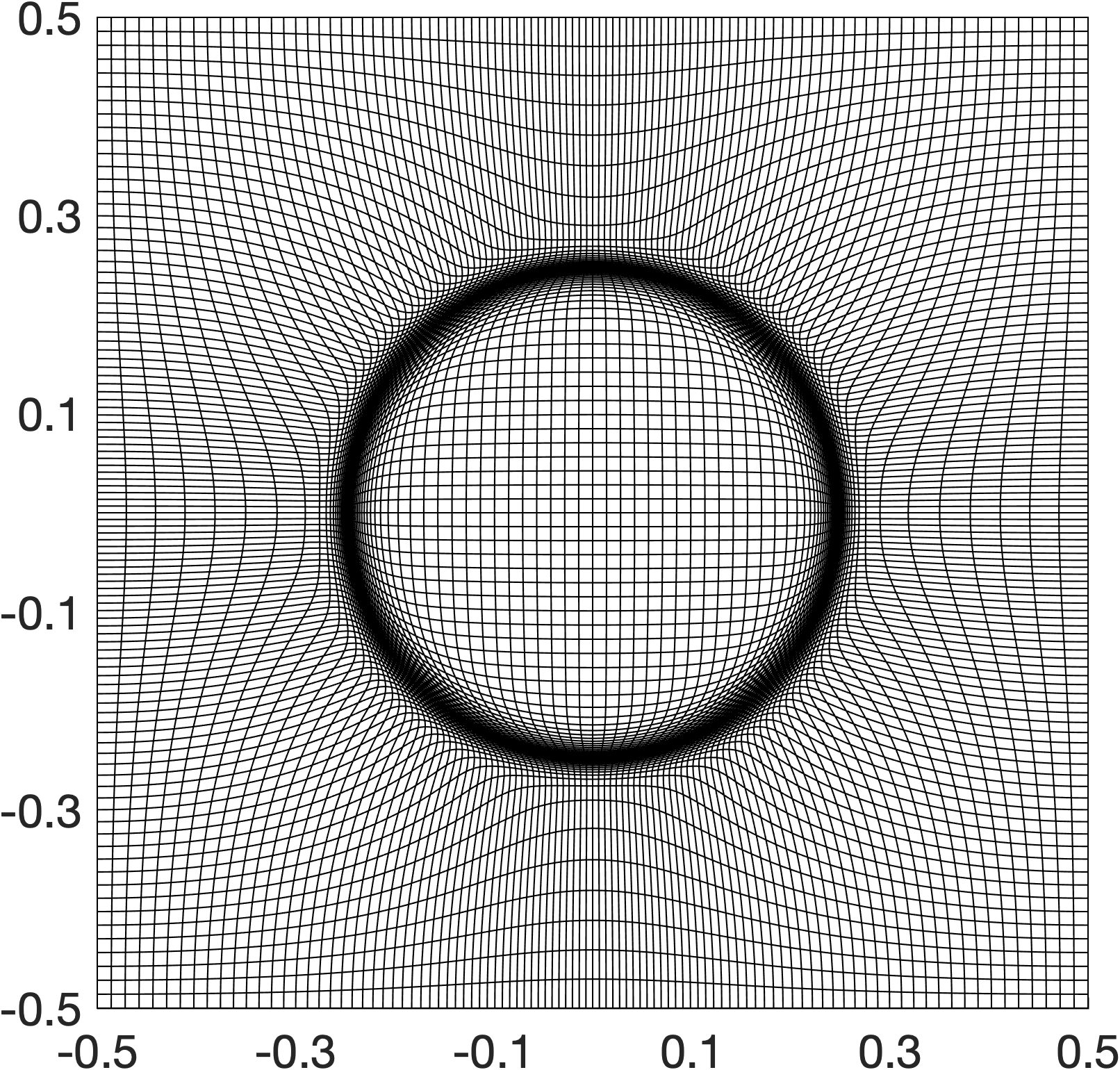}
		\caption{$\bm \alpha_3 = (20, 200, 0.25)$}
	\end{subfigure}
	\caption{Three high-order meshes are generated for the density functions shown in Figure \ref{ex3fig1} using a uniform background mesh of $100 \times 100$ quadrilaterals with $p=3$.}
	\label{ex3fig3}
\end{figure}

\begin{figure}[hthbp]
	\centering
 \includegraphics[width=0.49\textwidth]{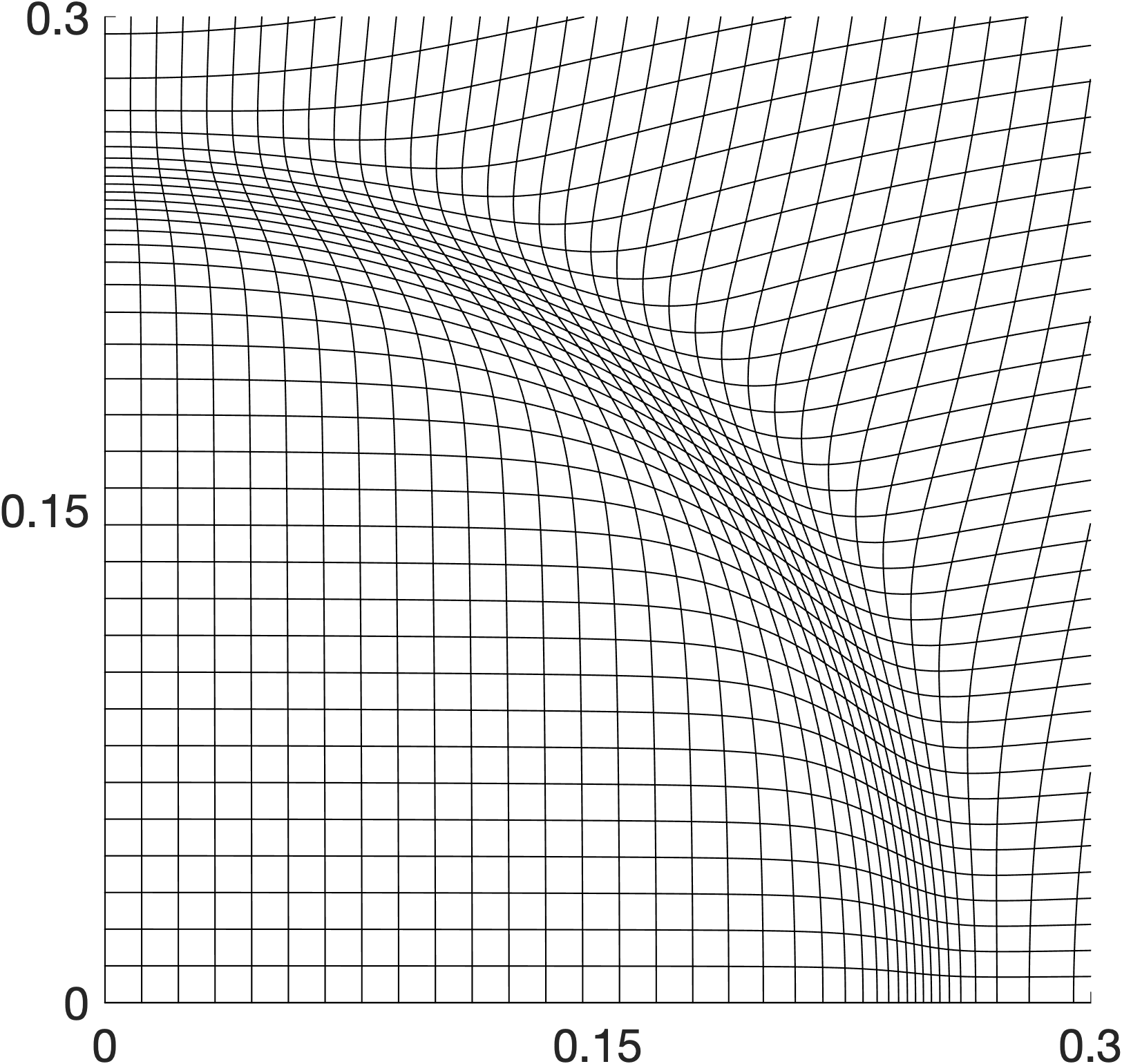} 
    \includegraphics[width=0.49\textwidth]{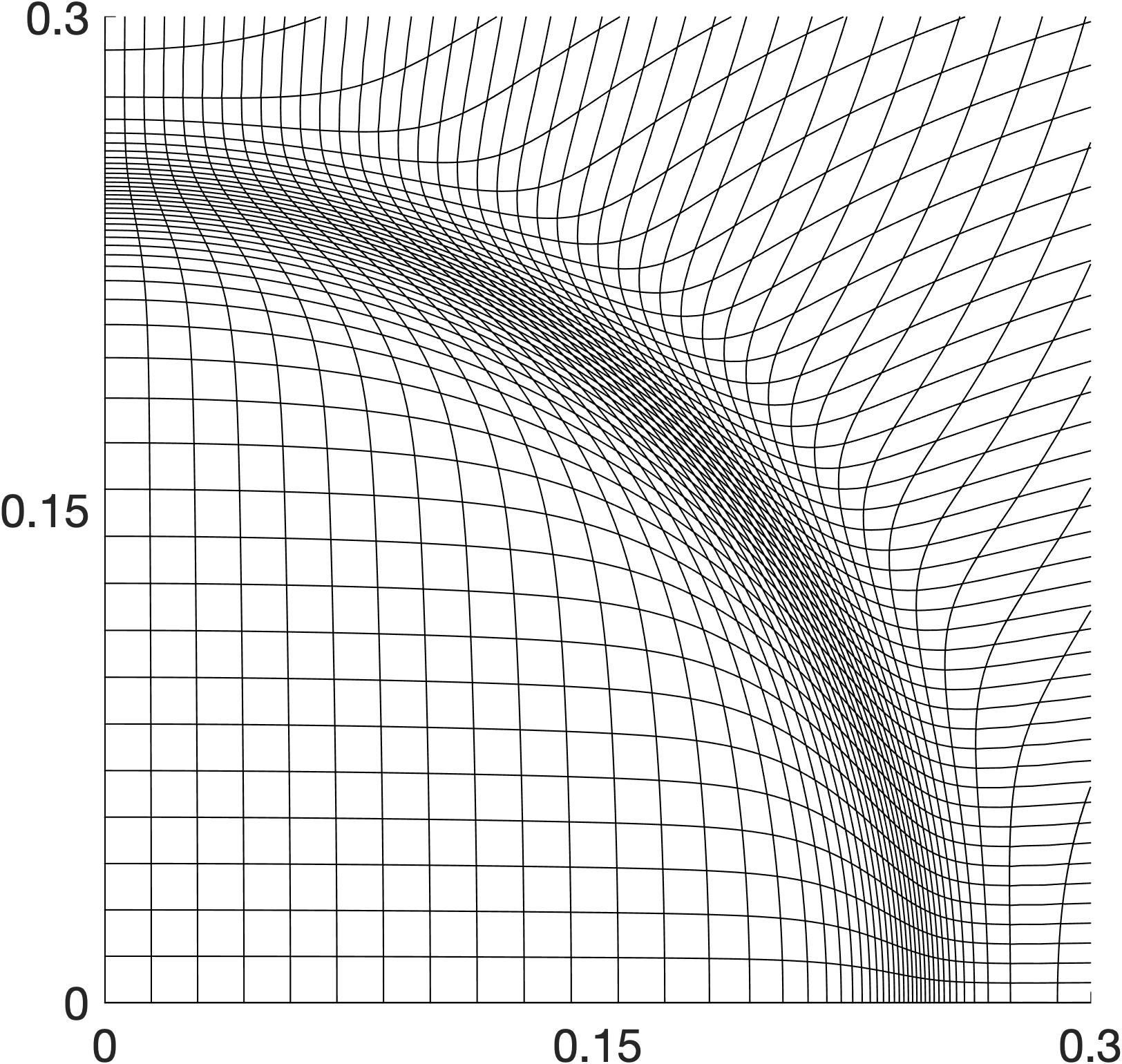}
	\caption{Close-up view near the ring of the meshes shown in Figure \ref{ex3fig3}(a) (left) and in Figure \ref{ex3fig3}(c) (right).}
	\label{ex3fig4}
\end{figure}

To demonstrate that our methods can generate smooth high-order meshes even at very high mesh resolutions, we consider a background mesh as a uniform grid of $100 \times 100$ quadrilaterals with polynomial degree $p=3$. Figure \ref{ex3fig3} depicts the three high-order meshes generated on this background mesh, while Figure \ref{ex3fig4} shows the close-up view near the ring of the first and third meshes. Despite there are many elements concentrating into the ring, the meshes are smooth and non-tangled.

\

\begin{figure}[hthbp]
	\centering
	\begin{subfigure}[b]{0.32\textwidth}
		\centering		\includegraphics[width=\textwidth]{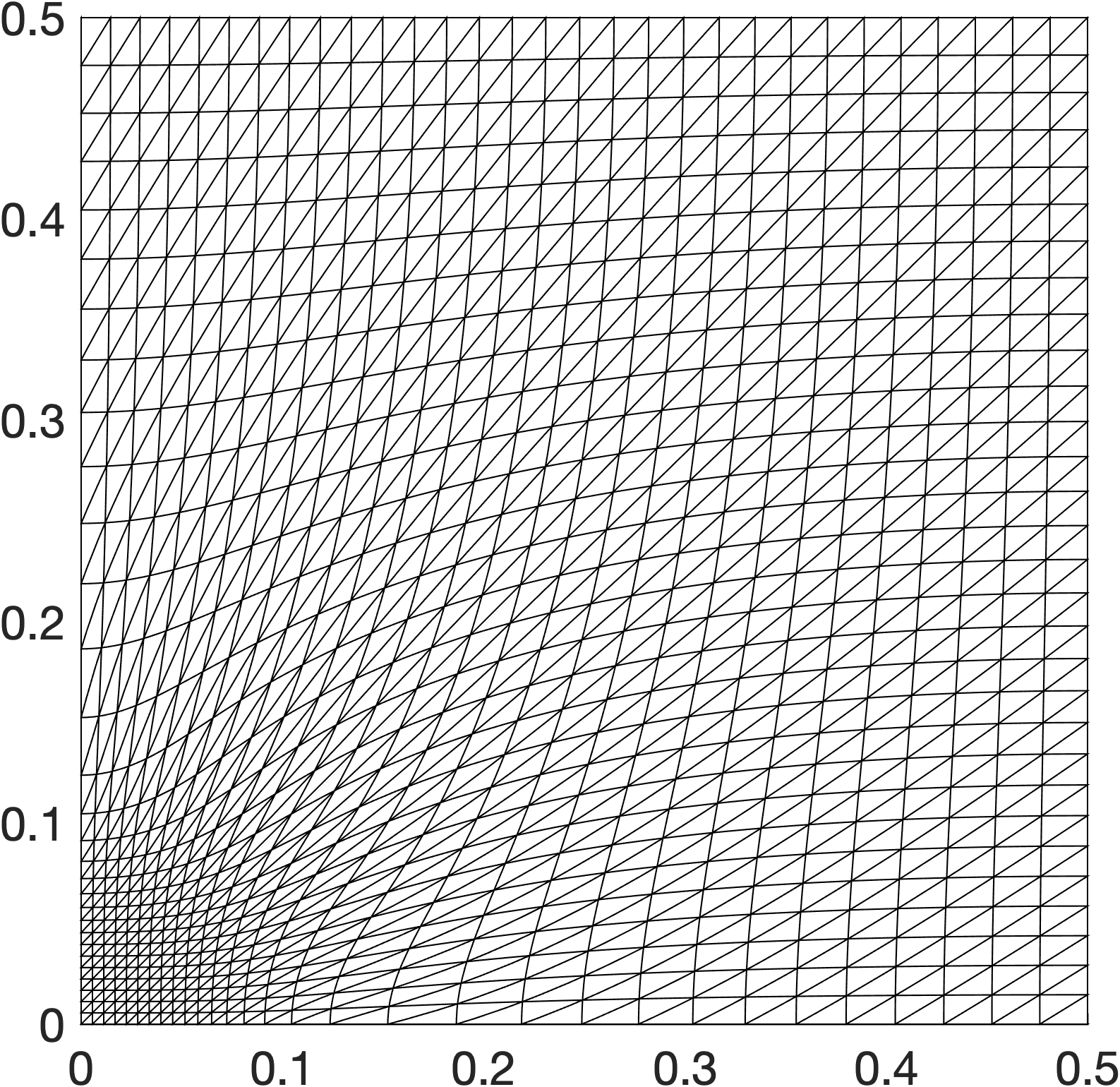}
		\caption{$\bm \alpha_4 = (10, 200, 0)$}
	\end{subfigure}
	\hfill
	\begin{subfigure}[b]{0.32\textwidth}
		\centering		\includegraphics[width=\textwidth]{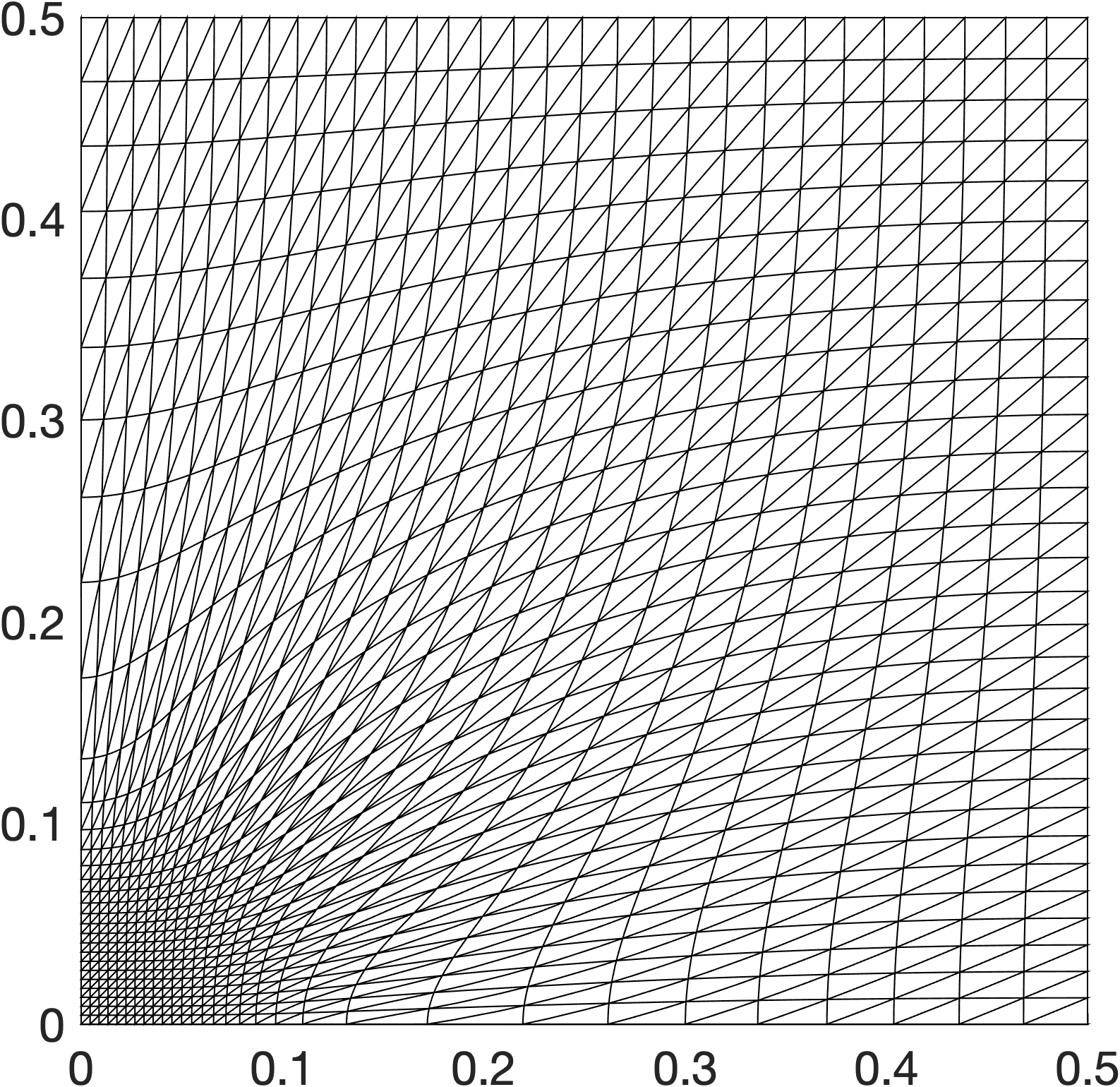}
		\caption{$\bm \alpha_5 = (20, 200, 0)$}
	\end{subfigure}
        \hfill
	\begin{subfigure}[b]{0.32\textwidth}
		\centering		\includegraphics[width=\textwidth]{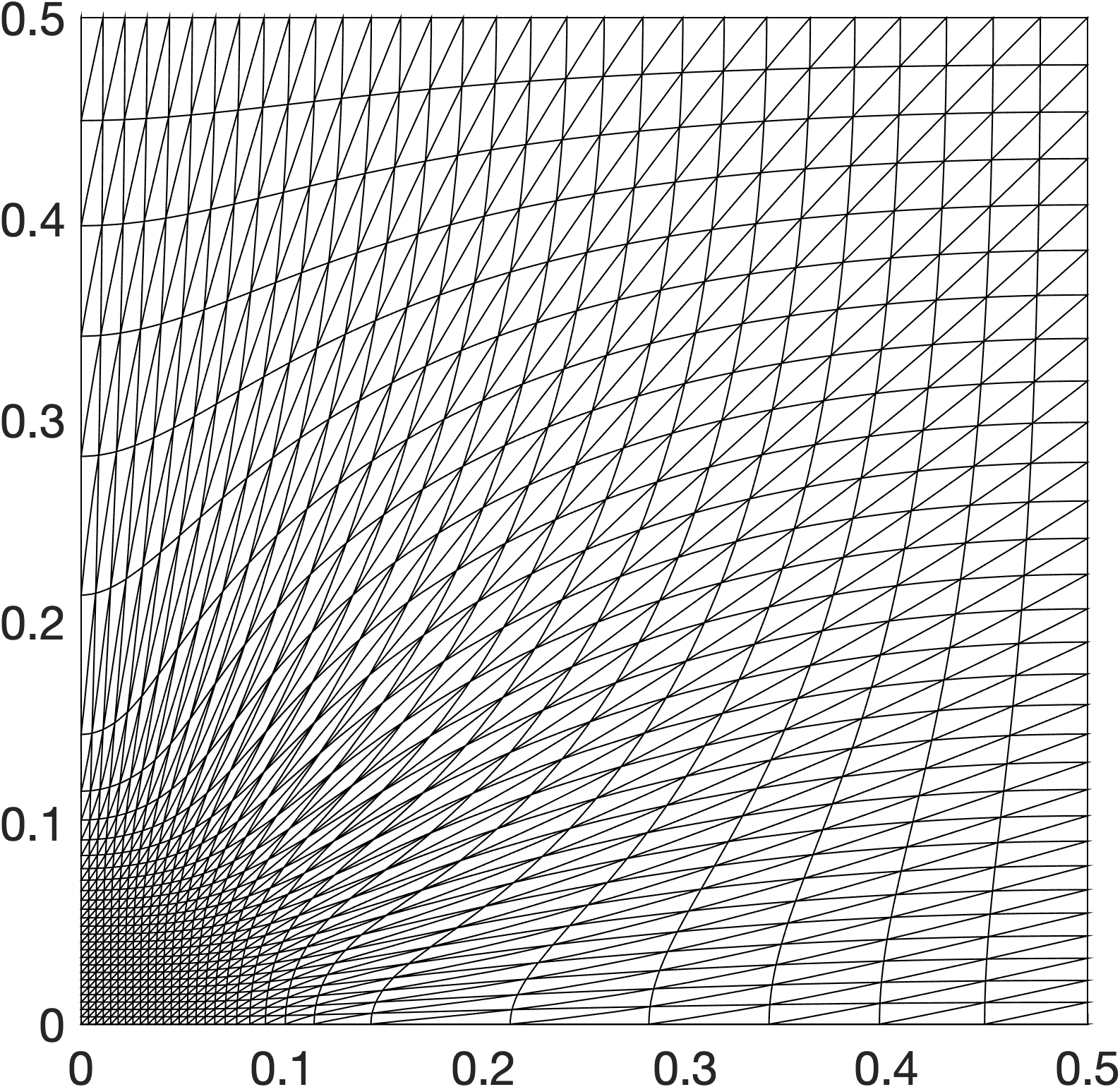}
		\caption{$\bm \alpha_6 = (40, 200, 0)$}
	\end{subfigure}
	\caption{Three high-order meshes are generated for three instances of the density function (\ref{dens}) using a uniform background mesh of $60 \times 60 \times 2$ triangles with $p=3$.}
	\label{ex4fig2}
\end{figure}

\begin{figure}[hthbp]
	\centering
	\begin{subfigure}[b]{0.32\textwidth}
		\centering		\includegraphics[width=\textwidth]{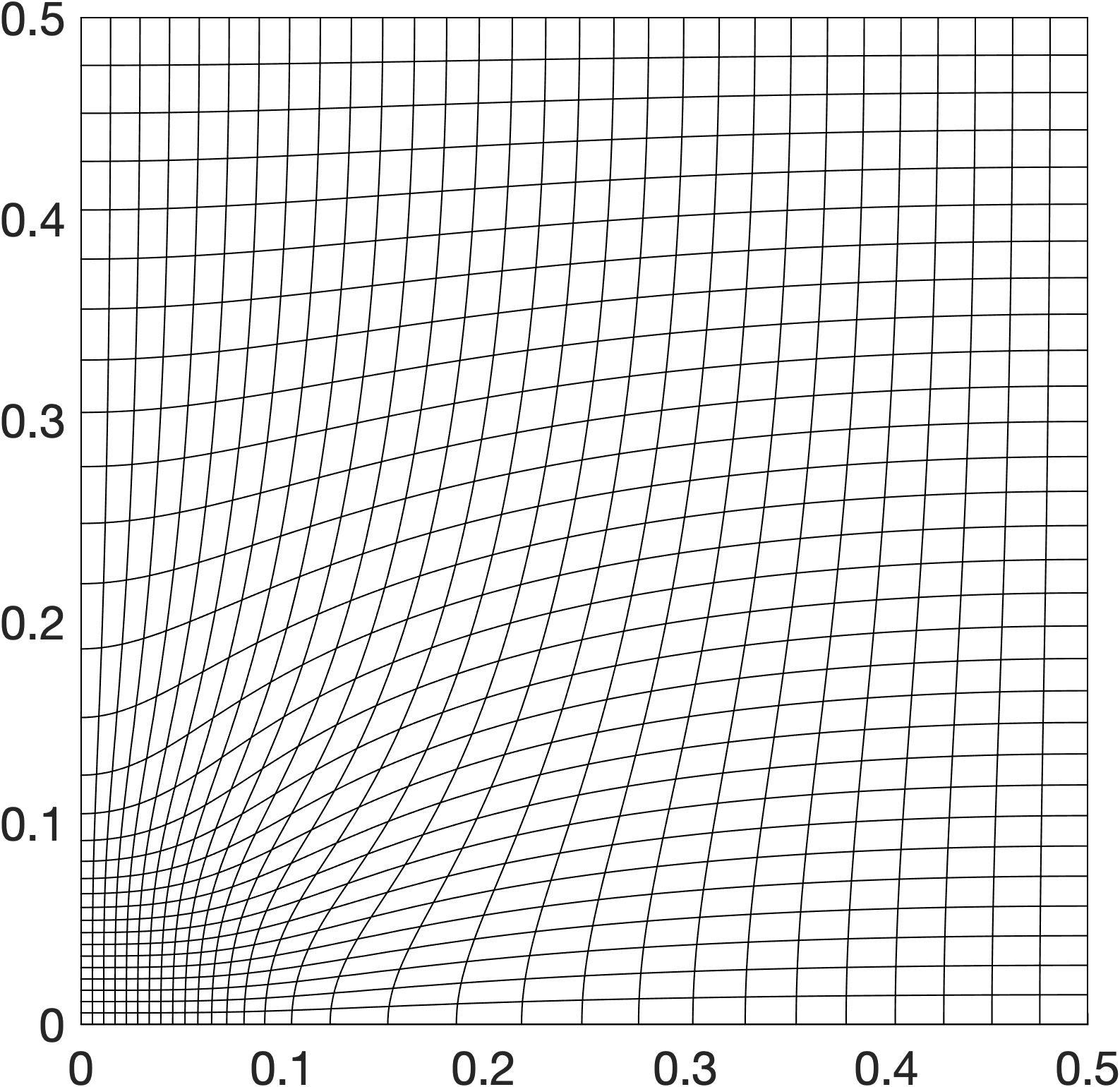}
		\caption{$\bm \alpha_4 = (10, 200, 0)$}
	\end{subfigure}
	\hfill
	\begin{subfigure}[b]{0.32\textwidth}
		\centering		\includegraphics[width=\textwidth]{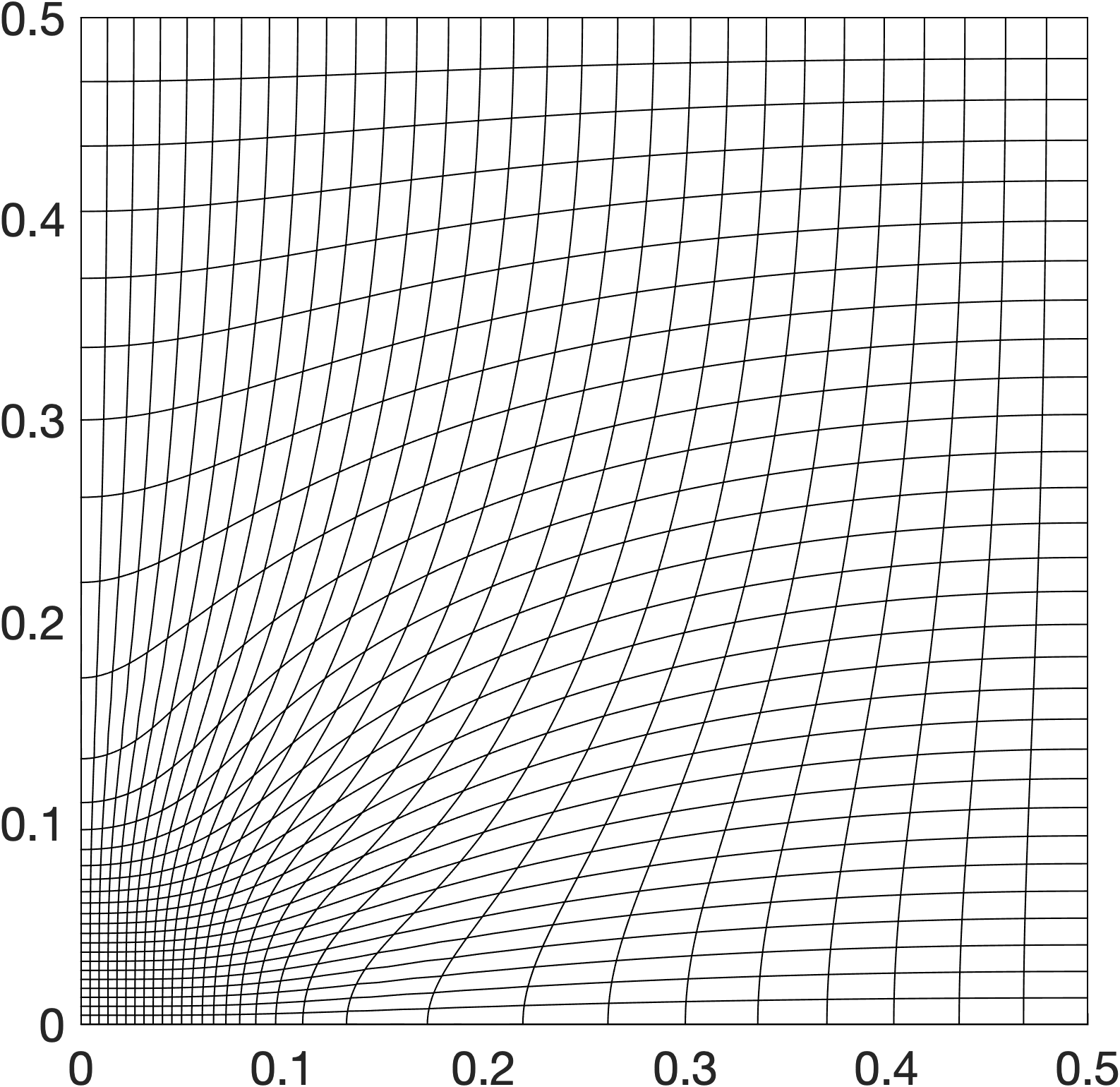}
		\caption{$\bm \alpha_5 = (20, 200, 0)$}
	\end{subfigure}
        \hfill
	\begin{subfigure}[b]{0.32\textwidth}
		\centering		\includegraphics[width=\textwidth]{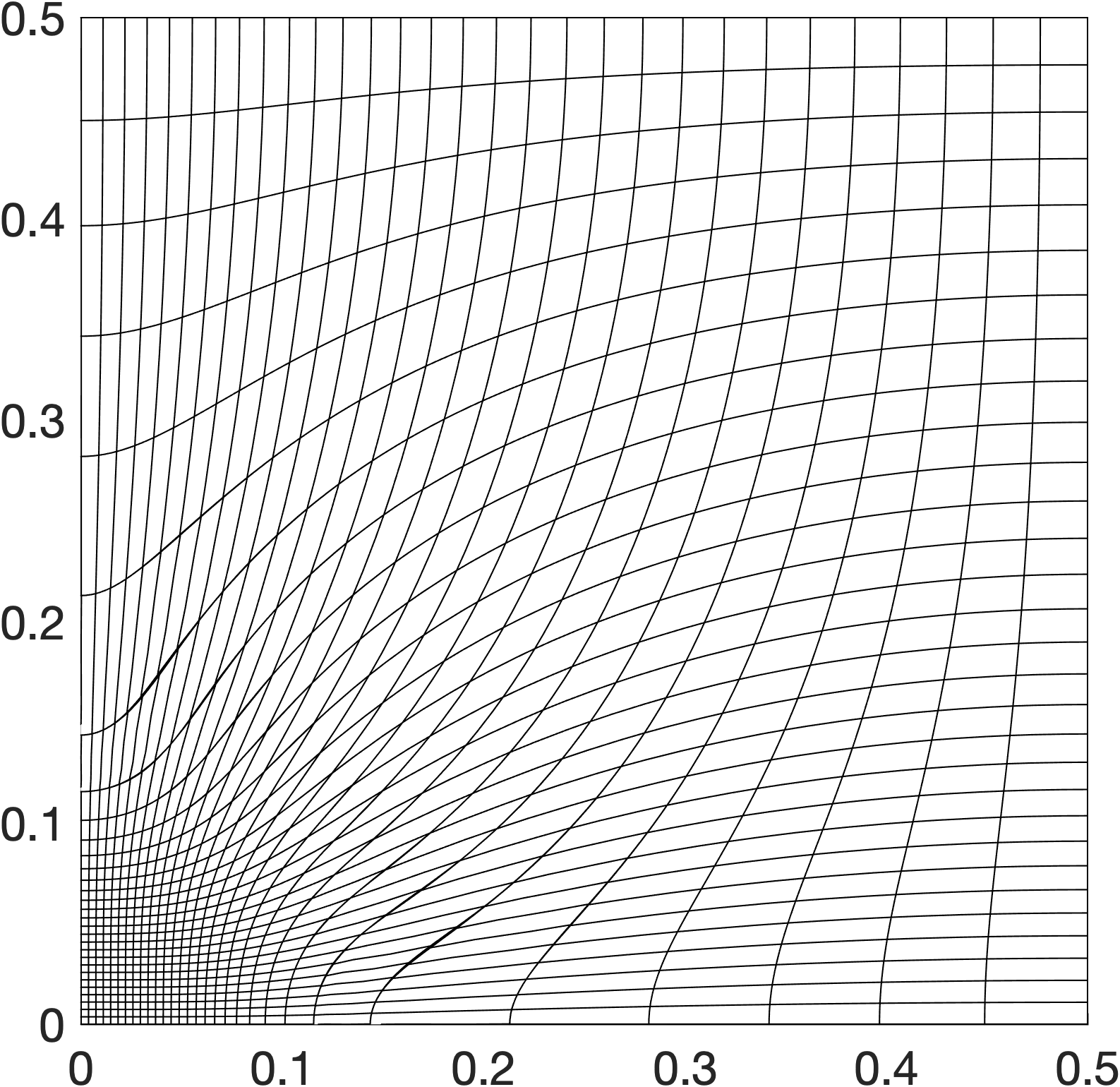}
		\caption{$\bm \alpha_6 = (40, 200, 0)$}
	\end{subfigure}
	\caption{Three high-order meshes are generated for three instances of the density function (\ref{dens}) using a uniform background mesh of $60 \times 60$ quadrilaterals  with $p=3$.}
	\label{ex4fig3}
\end{figure}

\noindent 
\textbf{Bell meshes on a square domain.}  We consider three new instances of the density function (\ref{dens}) corresponding to $\bm \alpha_4 = (10, 200, 0)$, $\bm \alpha_5 = (20, 200, 0)$, and $\bm \alpha_6 = (40, 200, 0)$. Figure \ref{ex4fig2} depicts the three high-order meshes that are generated for these instances on a uniform mesh of $60 \times 60 \times 2$ triangles with $p=3$, while Figure \ref{ex4fig3} shows the meshes on a uniform mesh of $60 \times 60$ quadrilaterals. We see that as $\alpha_1$ increases from 10, 20, to 40, it results in more elements concentrating into the origin $(0,0)$.

\begin{table}[ht]
  \begin{center}
\scalebox{0.8}{%
    $\begin{array}{|c||c c || c c|}
    \hline    
    & \multicolumn{2}{|c||}{\mbox{Newton-HDG method}} & \multicolumn{2}{|c|}{\mbox{Fixed point-HDG method}} \\ 
   \hline 
     \mbox{Case} & \mbox{triangular grid}  & \mbox{quadrilateral grid} & \mbox{triangular grid}  & \mbox{quadrilateral grid} \\
  \hline    
\bm \alpha_1 & 6  &  6  &  40  &  55  \\  
\bm \alpha_2 & 9  &  9  &  57  &  92  \\  
\bm \alpha_3 & 16  &  11  &  79  &  138  \\   
\bm \alpha_4 & 8  &  7  &  65  &  65  \\  
\bm \alpha_5 & 9  &  10  &  88  &  87  \\  
\bm \alpha_6 & 12  &  15  &  126  &  122  \\  
\hline
     \end{array} $
}
\caption{Number of iterations required to reach convergence for the HDG methods for different instances of the density function (\ref{dens}) on uniform triangular and quadrilateral grids.}
\label{ex34tab1}
\end{center}{$\phantom{|}$}     
\end{table}

\

\begin{figure}[hthbp]
	\centering
	\begin{subfigure}[b]{0.19\textwidth}
		\centering		\includegraphics[width=\textwidth]{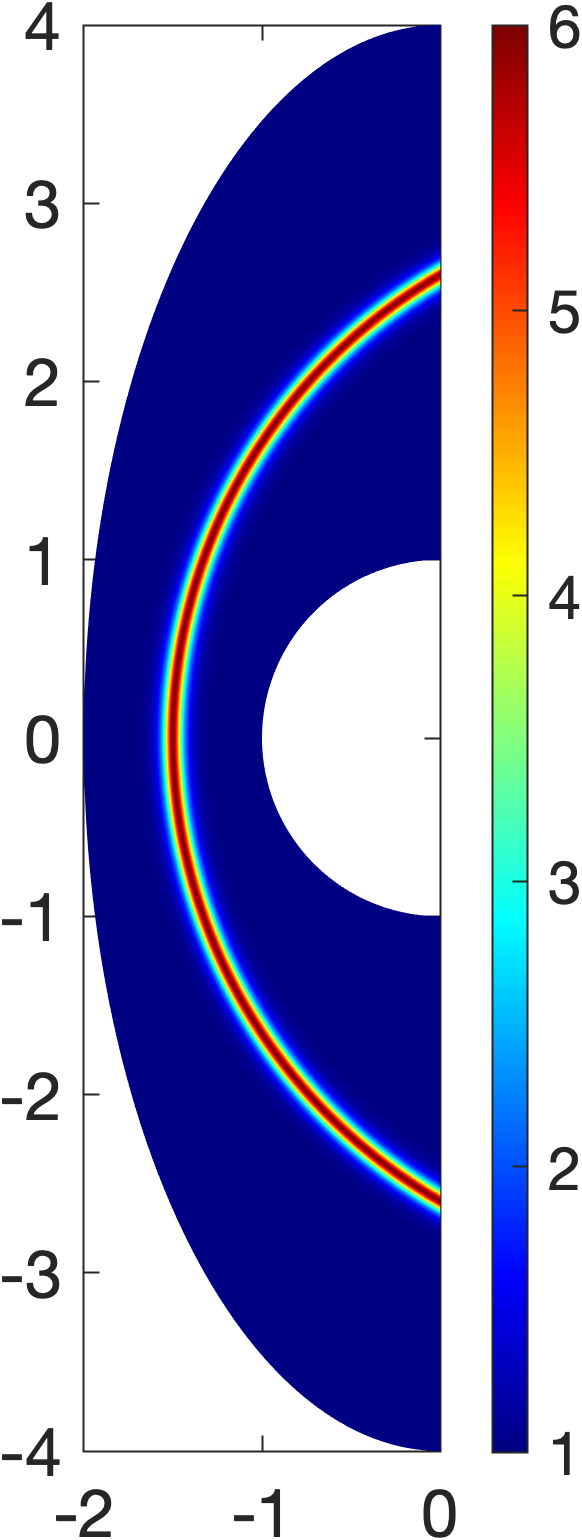}
		\caption{$\bm \beta_1 = (5, 5, 3)$}
	\end{subfigure}
	\hfill
	\begin{subfigure}[b]{0.19\textwidth}
		\centering		\includegraphics[width=\textwidth]{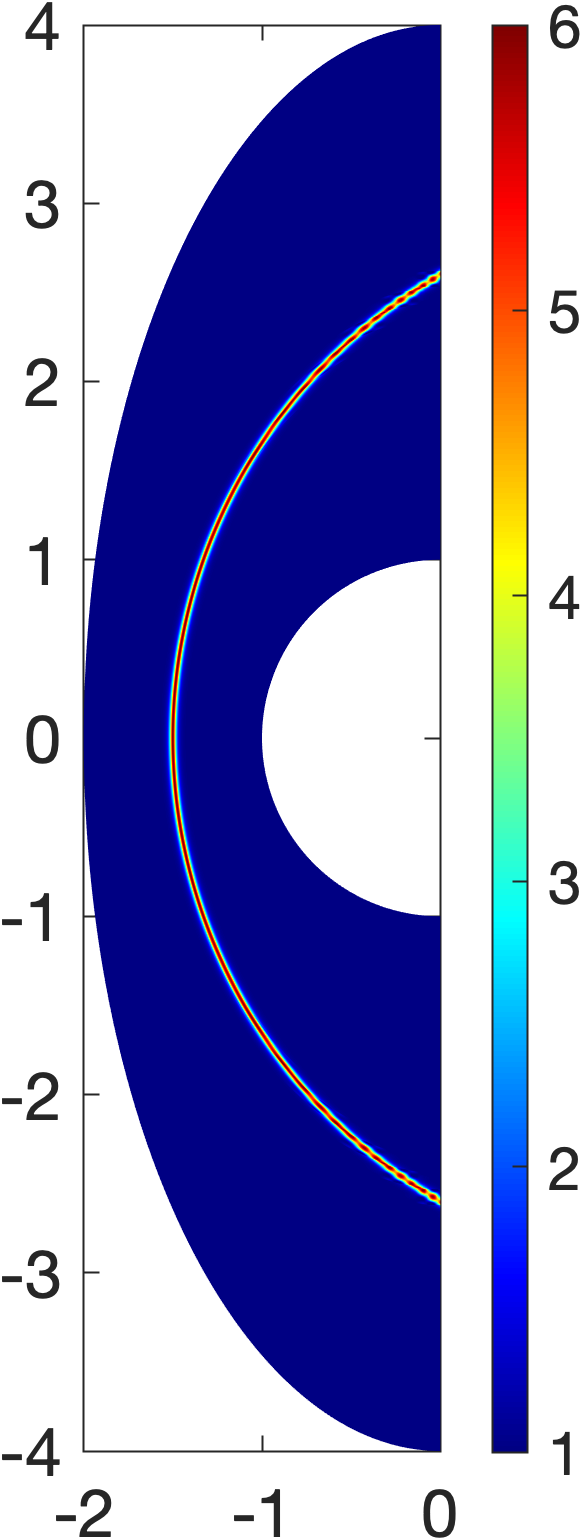}
		\caption{$\bm \beta_2 = (5, 15, 3)$}
	\end{subfigure}
        \hfill
	\begin{subfigure}[b]{0.2\textwidth}
		\centering		\includegraphics[width=\textwidth]{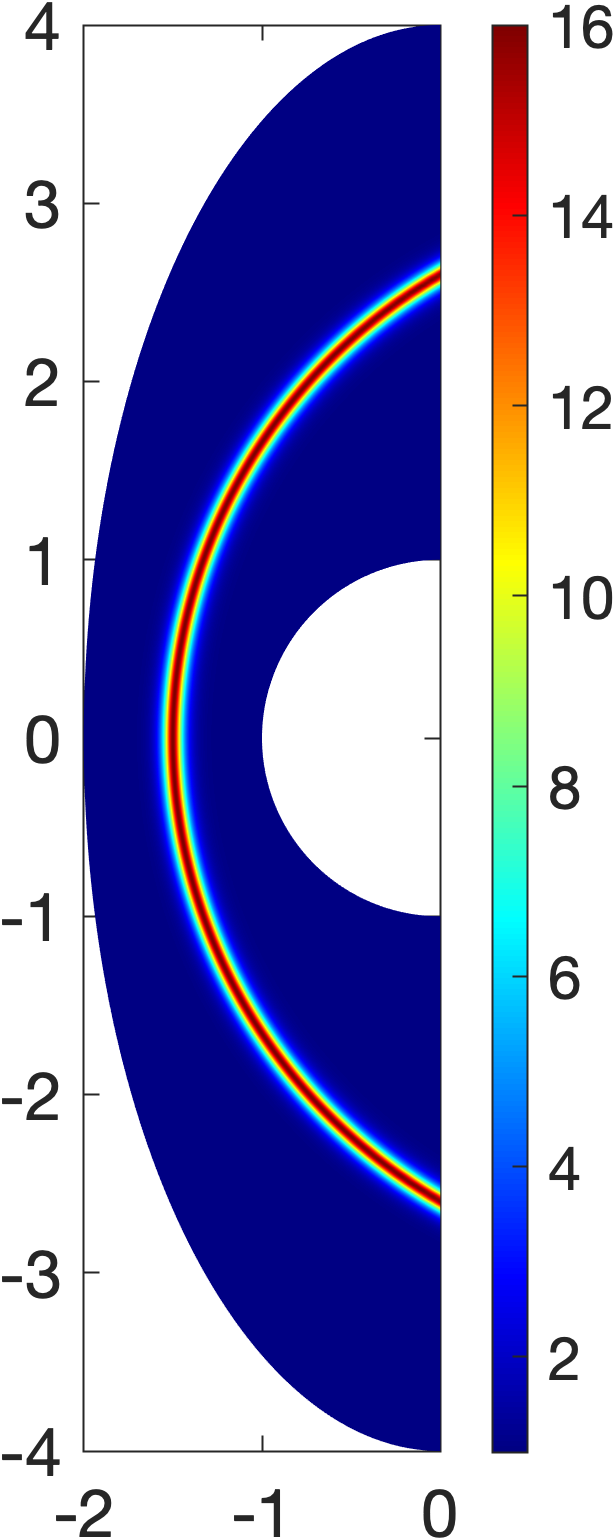}
		\caption{$\bm \beta_3 = (15, 5, 3)$}
	\end{subfigure}
        \hfill
	\begin{subfigure}[b]{0.2\textwidth}
		\centering		\includegraphics[width=\textwidth]{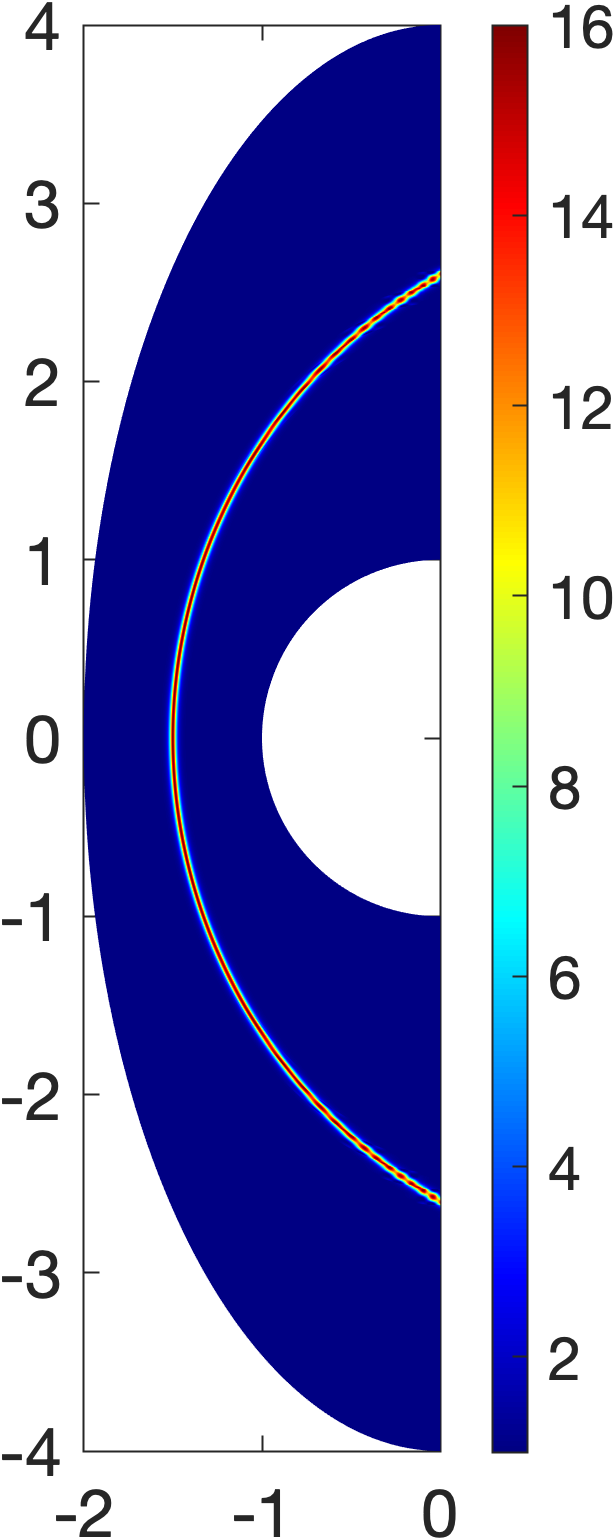}
		\caption{$\bm \beta_3 = (15, 15, 3)$}
	\end{subfigure}
	\caption{Four instances of the density function.}
	\label{ex5fig1}
\end{figure}

\noindent 
\textbf{Shock-aligned meshes on a cylindrical domain.} High speed flows past a unit circular cylinder is a popular test case in computational fluid dynamics \cite{Bai2022a,Barter2010,Ching2019,Nguyen2011a,Moro2016,persson06:_shock_capturing,Persson2013}. For high Mach numbers, a strong bow shock forms in front of the cylinder. Therefore, it is important to generate high-quality meshes to align the bow shock. The geometry is described by a half unit cylinder $x^2 + y^2 = 1$ and a half elliptical boundary $x^2/2^2 + y^2/4^2 = 1$. To represent a bow shock, the following target density function is considered 
\begin{equation}
\label{dens2}
\rho'(x', y') = 1 + \beta_1 \mbox{sech}^2 \left( \beta_2 \left( (x' - 1.5)^2 + y'^2 - \beta_3^2 \right) \right)  ,  
\end{equation}
where $\bm \beta = (\beta_1, \beta_2, \beta_3)$ determines the density function.  We consider four instances of the density function (\ref{dens2}) corresponding to $\bm \beta_1 = (5, 5, 3)$, $\bm \beta_2 = (5, 15, 3)$, $\bm \beta_3 = (15, 5, 3)$, and $\bm \beta_4 = (15, 15, 3)$, as shown in Figure \ref{ex5fig1}.

\begin{figure}[hthbp]
	\centering
	\begin{subfigure}[b]{0.14\textwidth}
		\centering		\includegraphics[width=\textwidth]{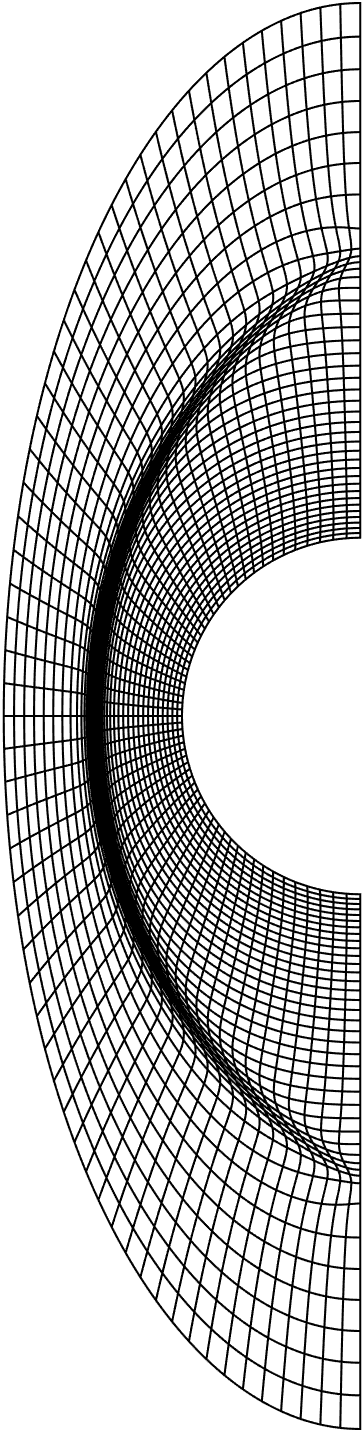}
		\caption{$\bm \beta_1 = (5, 5, 3)$}
	\end{subfigure}
	\hfill
	\begin{subfigure}[b]{0.14\textwidth}
		\centering		\includegraphics[width=\textwidth]{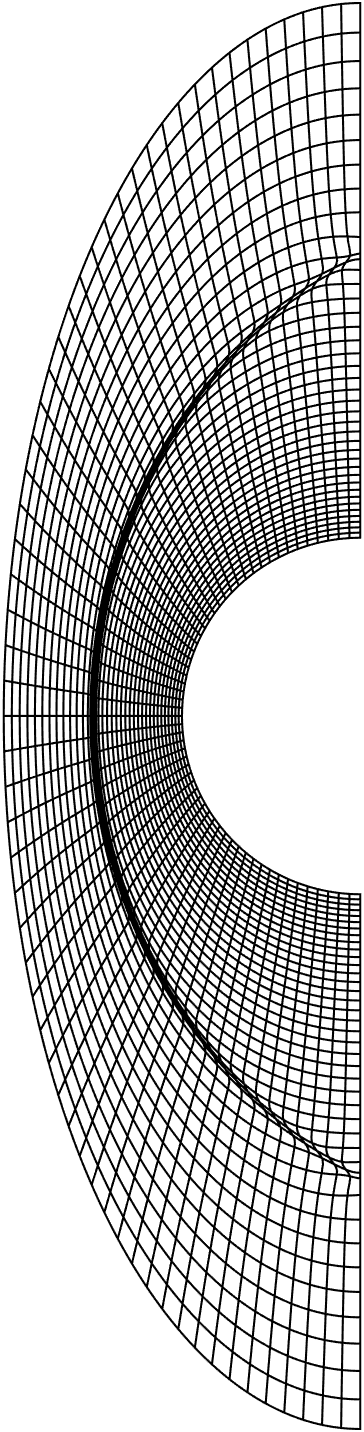}
		\caption{$\bm \beta_2 = (5, 15, 3)$}
	\end{subfigure}
        \hfill
	\begin{subfigure}[b]{0.14\textwidth}
		\centering		\includegraphics[width=\textwidth]{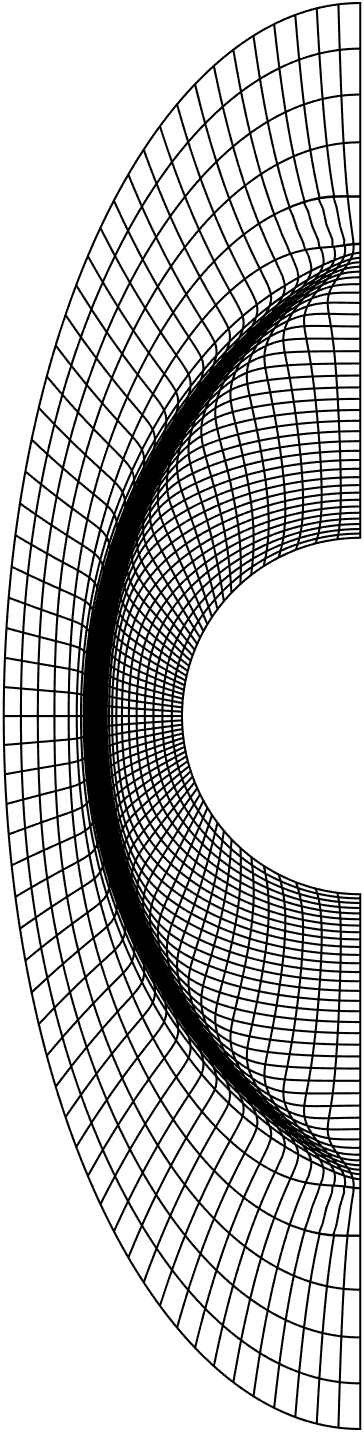}
		\caption{$\bm \beta_3 = (15, 5, 3)$}
	\end{subfigure}
        \hfill
	\begin{subfigure}[b]{0.15\textwidth}
		\centering		\includegraphics[width=0.935\textwidth]{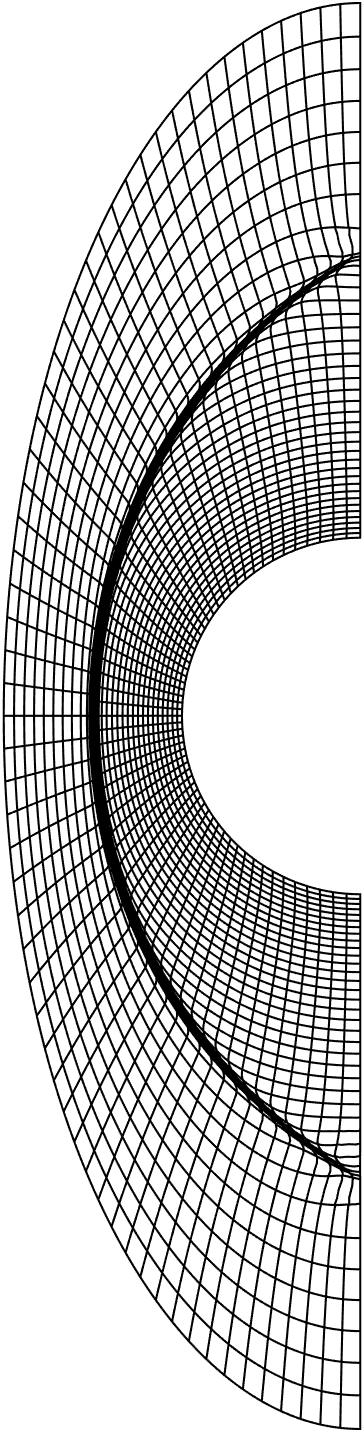}
		\caption{$\bm \beta_3 = (15, 15, 3)$}
	\end{subfigure}
	\caption{High-order quadrilateral meshes correspond to the four density instances shown in Figure \ref{ex5fig1}.}
	\label{ex5fig2}
\end{figure}

\begin{figure}[hthbp]
	\centering
 \includegraphics[width=0.47\textwidth]{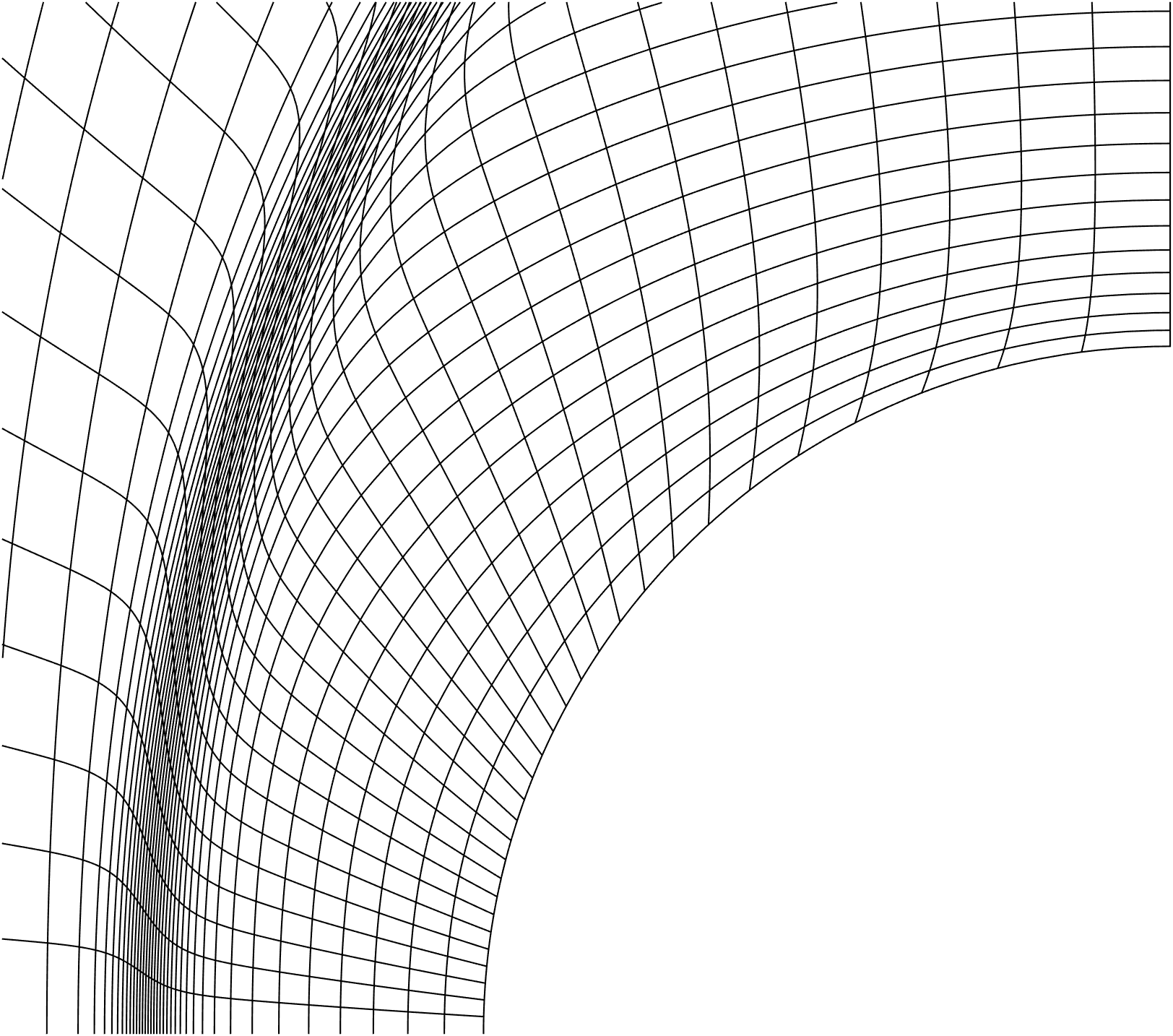} 
 \hfill 
    \includegraphics[width=0.47\textwidth]{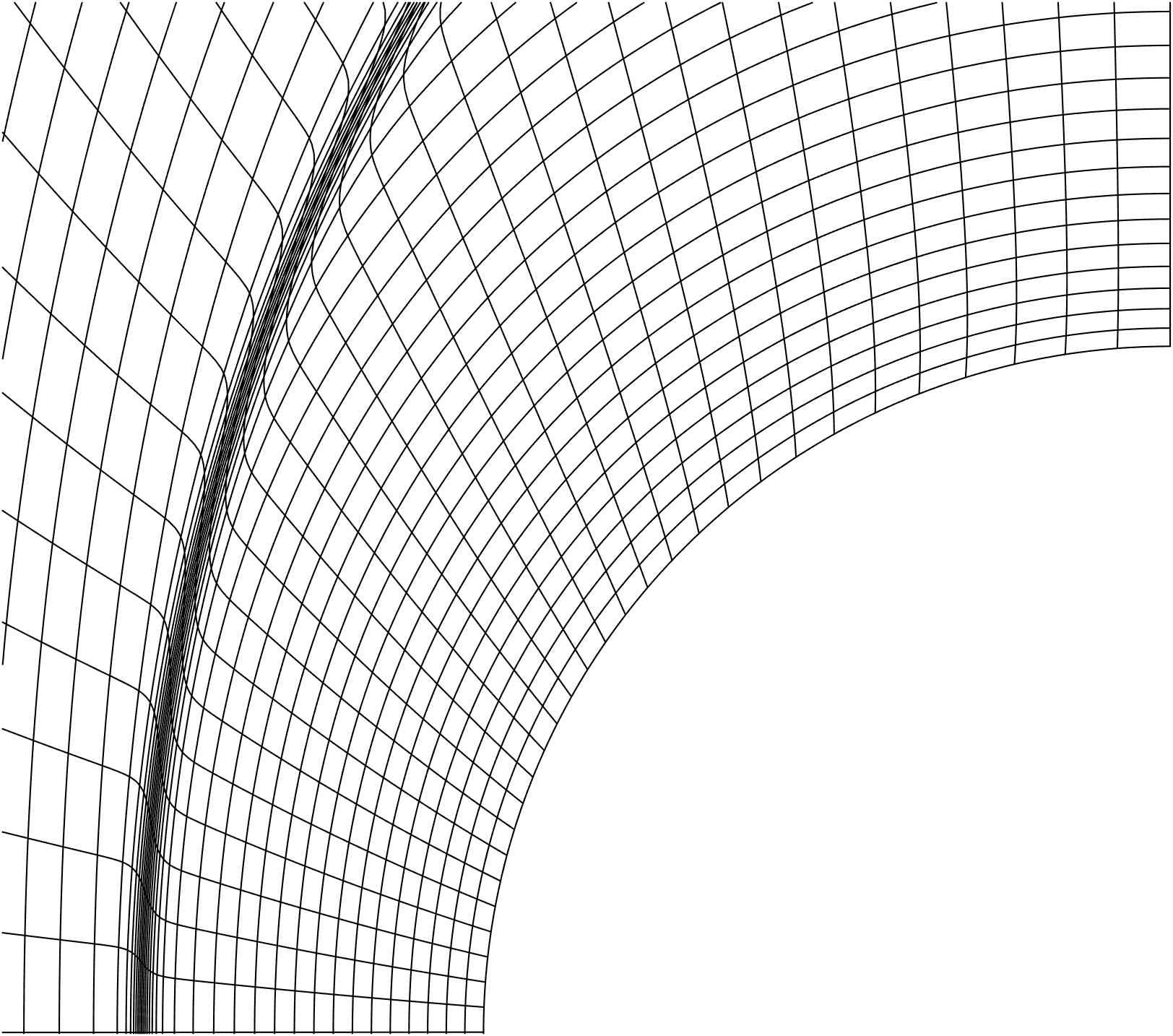}
	\caption{Close-up view near the cylinder of the meshes shown in Figure \ref{ex5fig2}(c) (left) and in Figure \ref{ex5fig2}(d) (right).}
	\label{ex5fig3}
\end{figure}

The background mesh on which the optimal transport meshes are generated is a quadrilateral grid of $40 \times 60$ elements with polynomial degree $p=4$. Figure \ref{ex5fig2} depicts the four high-order meshes generated for these density instances shown in Figure \ref{ex5fig1}, while Figure \ref{ex5fig3} shows the close-up view near the shock region of the last two meshes. We see that increasing the amplitude of the density function in the shock region results in more elements concentrating into the shock region. Furthermore, widening the density function increases the thickness of the shock region. We emphasize that the generated meshes are high-order, smooth, non-tangled, and conforming to the curved boundary of the physical domain. 

\section{Conclusion}



We have presented two hybridizable discontinuous Galerkin methods for numerically solving the Monge–Amp\`ere equation. The first HDG method is based on the Newton method, whereas the second HDG method is based on the fixed-point method. Numerical results were presented to demonstrate that the convergence and accuracy of the HDG methods. The Newton-HDG method is more efficient since it requires less iterations than the fixed-point HDG method. The numerical results showed that the HDG methods yield the convergence rate of $O(h^p)$ for the approximate scalar variable $u_h$, the approximate gradient $\bm q_h$, and the approximate Hessian $\bm H_h$. The convergence rate of $O(h^p)$ for the approximate gradient and Hessian is an attractive feature of the HDG methods. Furthermore, the HDG methods were extended to generate $r$-adaptive high-order meshes based on equidistribution of a given scalar density function via the optimal transport theory. Several numerical experiments were presented to illustrate the generation of smooth high-order meshes on planar and curved domains.  

 
It is important to analyze the HDG methods to understand the convergence rates observed in this paper. Furthermore, we would like to extend our methodology to generate $r$-adaptive meshes for flow problems in computational fluid dynamics. Therefore, three-dimensional mesh generation based on the optimal transport theory is an important topic to be addressed in future work.


\section*{Acknowledgements} \label{}

We would like to thank Professor Bernardo Cockburn at the University of Minnesota for fruitful discussions.  We gratefully acknowledge the United States  Department of Energy under contract DE-NA0003965, the National Science Foundation for supporting this work (under grant number NSF-PHY-2028125), and the Air Force Office of Scientific Research under Grant No. FA9550-22-1-0356 for supporting this work.  

 \bibliographystyle{elsarticle-num} 
\bibliography{library.bib}




\end{document}